\documentclass[final,leqno,onefignum,onetabnum]{siamltex1213}
\usepackage{amssymb}
\usepackage{amsmath}
\usepackage{theorem}
\usepackage{algorithm}
\usepackage{graphicx}
\usepackage{multirow}
\definecolor{red}{rgb}{.840,.352,.125}

\usepackage{soul,color}
\newif\ifshowamends
\newif\ifshownotes
\showamendstrue
\definecolor{delcol}{rgb}{0.9,0,0}
\definecolor{addcol}{rgb}{0,0.6,0}
\definecolor{notecol}{rgb}{0,0,0.9}
\definecolor{delcolhl}{rgb}{1,0.3,0}
\definecolor{addcolhl}{rgb}{0,1,0}
\definecolor{notecolhl}{rgb}{0.45,0.65,1}

\newcommand\NLadd[2][s]
	{\ifshowamends\ifx#1l{\color{addcol}#2}\else\sethlcolor{addcolhl}\hl{#2}\fi\else{#2}\fi}
\newcommand\NLdel[2][s]
	{\ifshowamends\ifx#1l{\color{delcol}#2}\else\sethlcolor{delcolhl}\hl{#2}\fi\fi}
\newcommand\nota[2][s]{\ifshownotes\ifx#1l{\begin{quote}\color{notecol}\slshape[#2]\end{quote}}\else\sethlcolor{notecolhl}\hl{\slshape[#2]}\fi\fi}

\newcommand{\R}{\mathbb{R}}

\newcommand{\Z}{\mathbb{Z}}
\newcommand{\ve}[1]{ #1}
\newcommand{\Tr}{\mbox{Tr}}
\newtheorem{Thm}{Theorem}

\newenvironment{ifelse}{%
  \begin{list}{}{%
      \setlength{\topsep}{0pt}\setlength{\parskip}{0pt}
      \setlength{\partopsep}{0pt}\setlength{\itemsep}{0pt}}}%
  {\end{list}}

\theoremheaderfont{\scshape}
\makeatletter

  \gdef\listctr{list\romannumeral\the\@listdepth}\expandafter
  
\makeatother

\newenvironment{AlgorithmSteps}[1][1]{%
  \begin{list}{\csname label\listctr\endcsname}{%
      \usecounter{\listctr}
      
      \settowidth{\labelwidth}{\textsc{Step\ #1.}}%
      \setlength{\leftmargin}{\labelwidth}\addtolength{\leftmargin}{\labelsep}}}%
  {\end{list}}

\def\xk{ x^{(k)}}
\def\xkk{ x^{(k+1)}}
\def\x{ x}

\def\aij{e}
\def\bij{g}

\title{A scaled gradient projection method for Bayesian learning in dynamical systems \thanks{This work has been partially supported by MIUR (Italian Ministry for University and Research), under the projects FIRB - Futuro in Ricerca 2012, contract RBFR12M3AC, and PRIN 2012, contract 2012MTE38N. The Italian GNCS - INdAM (Gruppo Nazionale per il Calcolo Scientifico - Istituto Nazionale di Alta Matematica) is also acknowledged.}}

\author{S. Bonettini\thanks{Dipartimento di Matematica e Informatica, Universit\`{a} di Ferrara, Via Saragat 1, 44121 Ferrara, Italy
(\email{silvia.bonettini@unife.it})} \and A. Chiuso\thanks{Dipartimento di Ingegneria dell'Informazione, Universit\`{a} di Padova, Via Gradenigo 6/b, 35131 Padova, Italy}\and M. Prato\thanks{Dipartimento di Scienze Fisiche, Informatiche e Matematiche, Universit\`{a} di Modena e Reggio Emilia, Via Campi 213/b, 41125 Modena, Italy}}

\begin{document}
\maketitle
\slugger{sisc}{xxxx}{xx}{x}{x--x}

\begin{abstract}
A crucial task in system identification problems is the selection of the most appropriate model class, and is classically addressed resorting to cross-validation or using order selection criteria based on  asymptotic arguments. As recently suggested in the literature, this can be addressed in a Bayesian framework, where model complexity is regulated by few hyperparameters, which can be estimated via marginal likelihood maximization. It is thus of primary importance to design effective optimization methods to solve the corresponding optimization problem. If the unknown impulse response is modeled as a Gaussian process with a suitable kernel, the maximization of the marginal likelihood leads to a challenging nonconvex optimization problem, which requires a stable and effective solution strategy.\\
In this paper we address this problem by means of a scaled gradient projection algorithm, in which the scaling matrix and the steplength parameter play a crucial role to provide a meaningful solution in a computational time comparable with second order methods. In particular, we propose both a generalization of the split gradient approach to design the scaling matrix in the presence of box constraints, and an effective implementation of the gradient and objective function.\\
The extensive numerical experiments carried out on several test problems show that our method is very effective in providing in few tenths of a second solutions of the problems with accuracy comparable with state-of-the-art approaches. Moreover, the flexibility of the proposed strategy makes it easily adaptable to a wider range of problems arising in different areas of machine learning, signal processing and system identification.
\end{abstract}



\pagestyle{myheadings}
\thispagestyle{plain}
\markboth{S. Bonettini, A. Chiuso, M. Prato}{A scaled gradient projection method for Bayesian learning in dynamical systems}

\section{Introduction}

System identification is concerned with automatic dynamic model building from measured data. Under this unifying umbrella, this field spans a rather broad spectrum of topics, considering different model classes (linear, hybrid, nonlinear, continuous and discrete time) as well as a variety of methodologies and algorithms, bringing together in a nontrivial way concepts from classical statistics, machine learning and dynamical systems.

The demand for reliable automatic tools for data based modeling of dynamical systems has attracted a considerable interest in the automatic control as well as in the statistics and econometrics communities since the $60'$s and has been mainly developed following the parametric maximum likelihood (ML)/prediction error (PE) framework, whose widespread use is to be attributed mainly to its attractive asymptotic statistical properties \cite{Ljung:99,Soderstrom,BoxJenkins}. Even if we restrict to linear, time-invariant, finite ``order'' dynamical systems (i.e. systems described by linear differential or difference equations with constant coefficients), where parametric methods are by now well developed and understood (see \cite{Ljung:99,Soderstrom}), it is fair to say that modeling cannot still be considered a ``completely automated'' task. For instance, in advanced process control applications \cite{YucaiZhu_Book}, modeling still is, by far, the most time consuming and costly step \cite{LHO_EJC2011}. As such, the demand for fast and reliable automated procedures for system identification makes this exciting field still a very active and lively one.

The system identification community, inspired by work in statistics \cite{Lasso1996,McKayARD}, machine learning \cite{Rasmussen,Tipping2001,Bach_MKL_2004} and signal processing \cite{Donoho2006,Wipf_IEEE_TIT_2011}, has recently developed and adapted methods based on regularization to jointly perform model selection and estimation in a computationally efficient and statistically robust manner \cite{SS2010,SS2011,ChenOL12,ChiusoPAuto2012,jmlr2014,SurveyKBsysid,AyazogluS12,FazelACC2001}. The main task of regularization is to control model complexity to face the so-called \emph{bias/variance dilemma} \cite{Ljung:99}. Different regularization strategies have been employed which can be classified in two main classes: regularization induced by smoothness priors (aka Tikhonov regularization, see \cite{KitagawaJASA1984,DoanLSER1984} for early references in the field of dynamical systems) and regularization for selection. This latter is usually achieved by convex relaxation of the $\ell_0$ quasi-norm (such as $\ell_1$ norm and variations thereof such as sum-of-norms, nuclear norm etc.) or other nonconvex sparsity inducing penalties which can be conveniently derived in a Bayesian framework, aka sparse Bayesian learning (SBL) \cite{McKayARD,Tipping2001,Wipf_IEEE_TIT_2011}.

In this paper we shall be concerned with regularization induced by smoothness priors; the structure of the chosen prior will bring in features usually encountered in SBL/automatic relevance determination \cite{McKayARD,Tipping2001} and multiple kernel learning \cite{Bach_MKL_2004,jmlr2014}. This makes the algorithms and results in this paper of a rather general interest. In particular we shall address the impulse response estimation problem for single input, single output (SISO) systems described by a convolution equation of the form
\begin{equation}\label{OEsystem}
y(t) = \sum_{k=1}^\infty h(k) u(t-k) + e(t),  \quad t \in \Z,
\end{equation}
where $y(t) \in \R$ is the \emph{output} signal, $u(t)\in \R$, is the measurable input signal, $h(k)$ is the (unknown) \emph{impulse response} and $e(t)$ is a zero mean white noise signal with unknown variance $\sigma^2$. As discussed in \cite{ChenetalTAC:14}, the very same framework studied in this paper can be easily adapted to identification of multi input single output (MISO) systems (see also \cite{ChiusoPAuto2012}), maintaining the key features which allow the application of the class of algorithms discussed herein.

We shall work in a Bayesian framework, thus modeling the unknown impulse response $h$ {(possibly an infinite dimensional object)} as a Gaussian process \cite{Rasmussen} with a suitable (prior) covariance $P(\nu)$ \cite{SS2010,SS2011,ChenOL12} (also known as \emph{kernel}). The chosen covariance is usually described by some unknown hyperparameters $\nu$ which give the prior enough flexibility to encode a sufficiently wide class of impulse responses.  The number of hyperparameters is typically small as compared to  the number of data as well as to the ``dimension'' of $h$, which as mentioned above can be infinite dimensional.  These hyperparameters can be estimated from data in a variety of ways; empirical evidence as well as some theoretical results \cite{AutomaticaML2014} support the use of the so-called \emph{marginal likelihood} (i.e. the data likelihood as a function of the unknown hyperparameters, having marginalized the unknown impulse responses from the joint density of data and unknowns) for hyperparameter estimation.
This boils down to a challenging optimization problem with the following features:
\begin{itemize}
\item it is nonconvex;
\item it requires to handle a large number of data $(y,u)$ (also several thousands) even when the number of unknowns (hyperparameters) is not too large (some tenths in most cases);
\item the Hessian matrix is, in some cases, quite costly to compute;
\item the computation of the objective function and its gradient requires the factorization of matrices which can be extremely ill-conditioned.
\end{itemize}
Thus, stable and effective algorithms should be designed carefully taking into account the features of the problem.
In particular, the simple structure of the constraints, which usually reduce to non-negativity or box, can be exploited by suitable projection methods.\\
In this paper we propose a scaled gradient projection method for marginal likelihood optimization, whose basic ingredients are the variable stepsize and scaling matrix, which are computed with a negligible computational cost at each gradient projection iteration. The stepsize parameter is chosen according to the Barzilai--Borwein rules, while the scaling matrix is based on a gradient decomposition technique. In spite of the theoretical convergence rate estimate, which in general classifies the classical gradient projection method as linearly convergent, it has been shown in the recent literature that the combination of these choices makes it a very practical, effective and robust numerical tool for several signal and image restoration problems \cite{BZZ09,Prato2012,Zanella2009,Zanella2013b}.\\
In this paper we show that, with a suitable choice of the scaling matrix and a careful implementation, the scaled gradient projection method applied to the impulse response estimation problem outperforms some second order state-of-the-art methods, leading to a significant reduction of the computational time.\\
The plan of the paper is the following. In Section \ref{sec:mod} we describe the system identification problem in the framework of the Bayesian approach, deriving the corresponding optimization problem, whose main features are described in Section \ref{sec:problem}. The proposed optimization method is presented in Section \ref{sec:opt}, focusing on steplength and scaling matrix selection. In particular, in Section \ref{sec:scaling}, we consider the split gradient strategy, which is a state-of-the-art approach for defining the scaling matrix in presence of non-negativity constraints, and we extend it to the more general case of box constraints. Some important implementation issues are discussed in Section \ref{sec:implementation}. Finally, the results of an extensive numerical experience are presented in Section \ref{sec:num}, showing the effectiveness of the proposed approach on the system identification problem, also with respect to other recent solvers. Our conclusions are offered in Section \ref{sec:conc}.

\subsection*{Notation} In the following the symbol $\mbox{Tr}(\cdot)$ indicates the matrix trace and  $\det(\cdot)$  the matrix determinant. We shall deal  with real random vectors whose (possibly conditional) measure, will always be absolutely continuous w.r.t. the Lebesgue measure and will thus admit a density $p$. We shall denote with $p(v)$ the density of $v$ (always w.r.t. the Lebesgue measure), $p(v|w)$ the conditional density of $v$ given $w$. Densities may depend upon some parameters (say $x$), in which case we shall use subscripts such as  $p_x(v)$ or $p_x(v|w)$.

\section{Problem statement and model derivation}\label{sec:mod}

We shall consider the following problem: given a finite data record $\{u(t),y(t)\}_{t=1}^N$ from system \eqref{OEsystem}, find an estimator of the impulse response $h$. This is clearly an ill-posed inverse problem since the unknown $h$ is an infinite dimensional object. As customary in the literature on inverse problems \cite{Bertero1} this can be tackled using Tikhonov regularization. Equivalently the (infinite dimensional) unknown $h$ can be modeled as a Gaussian process \cite{Rasmussen}. We shall follow this second route since it provides a natural way to introduce estimators of the \emph{regularization (hyper)parameters} through the marginal likelihood. We refer the reader to \cite{AutomaticaML2014} and references therein for some recent work in support of this approach.

{In order to {avoid theoretical issues related to dealing with infinite dimensional unknowns, chiefly the  complication of introducing probability densities for infinite dimensional objects, the  unknown impulse response} $\{h(k)\}_{k\in \Z}$ is truncated to a \emph{finite dimensional}, yet arbitrarily long vector. This approximation is always possible (within any arbitrary accuracy) since  the impulse response of a finite dimensional linear systems $\{h(k)\}_{k\in \Z^+}$ decays exponentially fast as a function of the index $k$. In addition, since only $N$ data points are available, no information could ever be obtained from data on the ``tail'' of the impulse response for $k\geq N$.
Thus the  model \eqref{OEsystem}
can be rewritten as
\begin{equation}\label{system1}
y(t) = \phi(t)^T\theta +e(t), \ \ \ t=n+1,\cdots,N \quad \theta\in\R^n
\end{equation}
where $\phi(t) = \left(u(t-1), u(t-2),\cdots,u(t-n)\right)^T$ and $\theta\in\R^n$ is the vector whose components are the system impulse response coefficients.\\ Note that, depending on the ``true'' underlying system, $n$ can  be arbitrarily large, so that estimating $\theta$ in the model \eqref{system1} is still an ill-conditioned inverse problem. We stress that this truncation is inessential; by resorting to \emph{Reproducing Kernel Hilbert Space} theory one can deal with the original infinite dimensional problem, see \cite{SS2010,SS2011}.}
Equation \eqref{system1} can be represented in matrix form as
\begin{equation}\label{LinModel}
Y = \Phi\theta + E,
\end{equation}
where $Y=(y(n+1),y(n+2),\ldots,y(N))^T$, $\Phi = (\phi(n+1)^T,\phi(n+2)^T,\ldots,\phi(N)^T)^T$ and $E =(e(n+1),e(n+2),\ldots,e(N))^T$. Since the noise affecting the data is white and Gaussian distributed with zero mean and variance $\sigma^2$, then  $Y$ conditioned on $\theta$ is Gaussian,  $Y|\theta \sim \mathcal{N}(\Phi\theta, \sigma^2I_{N-n})$, and thus has conditional density
\begin{equation*}
 \quad p_{\sigma^2}(Y|\theta)= (2\pi\sigma^2)^{-(N-n)/2}e^{-\frac{\|Y - \Phi\theta\|_2^2}{2\sigma^2}}.
\end{equation*}
We further model $\theta$ as a Gaussian random vector, independent of $E$, i.e.
\begin{equation*}
\theta \sim \mathcal{N}(\theta^{ap}, P(\nu)), \quad \quad p_\nu(\theta)=(2\pi)^{-n/2}\det(P(\nu))^{-1/2}e^{-\frac{1}{2}(\theta-\theta^{ap})^T P(\nu)^{-1}(\theta-\theta^{ap})},
\end{equation*}
where $P(\nu)$ is the prior covariance parametrized by the hyperparameter
vector $\nu\in\R^m$. Typical examples of prior covariance $P(\nu)$ will be given in Section \ref{sec:kernels}; suffices here to say that the number of hyper parameters $m$ is  typically ``small'' w.r.t. to the number of data points (from a few units to a few tens). For convenience of notation we shall define $\x=(\nu^T,\sigma^2)^T$. From \eqref{LinModel} it follows that $Y$ is the linear combination of independent Gaussian random vectors and, therefore,  the \emph{marginal likelihood} $p_\x(Y)$, {i.e. the marginal  of $Y$ obtained integrating $p_\x(Y,\theta)=p_{\sigma^2}(Y|\theta)p_{\nu}(\theta)$ w.r.t. $\theta$,} is still a multivariate normal  with mean $\Phi\theta^{ap}$ and covariance matrix
\begin{equation}\label{Sigma}
\Sigma(\x) = \Phi P(\nu)\Phi^T  + \sigma^2I_{N-n}, \quad \x=(\nu^T,\sigma^2)^T
\end{equation}
Using Bayes' Theorem we can compute the posterior density of $\theta$ given $Y$
\begin{equation}\label{BAYES}
p_{\x}(\theta | Y) = \frac{p_{\sigma^2}(Y|\theta)p_{\nu}( \theta)}{p_{\x}(Y)},
\end{equation}
which still depends on the unknown hyperparameters $x$. There are typically two approaches to deal with the unknown hyperparameters $x$. The first is the so called \emph{Full Bayes approach}: a prior distribution (possibly uninformative) for the hyperparameters is postulated which allows to integrate them out. The second,  which we consider in this paper, is the so-called  \emph{Empirical Bayes approach} \cite{Maritz:1989}: a point estimate $\widehat x$ of the hyperparameters $x$ is found and then the posterior \eqref{BAYES} is computed with $x$ fixed to its point estimate $\hat x$.
In this paper $\hat x$ is obtained following
 the maximum likelihood (ML) approach:
\begin{equation}\label{pr1}
\widehat{\x} = \underset{\x\in \Omega}{\rm{argmax}} \ p_\x(Y) = \underset{\x\in \Omega}{\rm{argmin}} \ f(\x),
\end{equation}
where $\Omega$ is some suitable subset of $\R^{m+1}$ and
\begin{align*}
f(\x) =& -2\log p_{\x}(Y) - (N-n)\log(2\pi) \\
=&  \log\det(\Sigma(\x)) + (Y-\Phi\theta^{ap})^T \Sigma(\x)^{-1}(Y-\Phi\theta^{ap}).
\end{align*}
After a solution $\widehat{\x}$ of problem \eqref{pr1} has been found, the maximum a posteriori (MAP) estimate $\widehat{\theta}$ of $\theta$, which is equal to the posterior mean for symmetric densities, can be computed:
%
\begin{align}
\widehat{\theta}: =& \ \underset{\theta}{\rm{argmax}} \ p_{\hat \x}(\theta|Y) = \  \underset{\theta}{\rm{argmin}} \ -2 \log(p_{\hat\x}(\theta|Y)) \nonumber \\
								 =& \ (\Phi \theta +\hat\sigma^2P(\hat\nu)^{-1})^{-1}(\Phi^T  Y + \hat\sigma^2P(\hat\nu)^{-1}\theta^{ap}),\label{thetahat}
\end{align}
where in the last equality the fact that $P(\hat\nu)$ is symmetric has been used.\\
Unless strong prior knowledge is available, the \emph{a priori} mean $\theta^{ap}$ is set to zero, thus estimating the impulse response coefficients $\theta$ requires going through the following steps:
\begin{enumerate}
\item solve the nonconvex, constrained optimization problem \eqref{pr1} where
\begin{equation}\label{fobj}
 f(\x) = Y^T\Sigma(\x)^{-1}Y+\log\det(\Sigma(\x))
\end{equation}
and $\Sigma(\x)\in \R^{(N-n)\times(N-n)}$ is defined in \eqref{Sigma}
\item compute the corresponding impulse response coefficients setting $\theta^{ap}=0$ in \eqref{thetahat}:
\begin{equation*}
\hat\theta = (\Phi \Phi^T + \hat\sigma^2P(\hat\nu)^{-1})^{-1}\Phi^T  Y\ .
\end{equation*}
\end{enumerate}
\subsection{Kernel matrices}\label{sec:kernels}
Several kernel matrices have been introduced in the recent years to model impulse responses of dynamical systems. Perhaps the major breakthrough has been the observation that the kernel has to capture structural properties of dynamical systems \cite{SS2010,DinuzzoKernels,ChiusoCLPCDC2014}, such as the fact that for linear systems described by difference/differential equations, the impulse response is a linear combination of exponentially decaying functions \cite{Kailath}. In order to do so, the seminal paper \cite{SS2010} has introduced the family of \emph{stable-spline} kernels; the most used kernels in this family are the \emph{stable-spline kernel of order 1}, called also tuned/correlated (TC) kernel
\cite{ChenOL12}:
\begin{subequations}\label{nonlinearkernels}
\begin{equation}\label{TCkernel0}P^{TC}_{k,j}(\nu) = c \cdot \min(\mu^{k},\mu^j),\ \ \ k,j=1,...,n, \end{equation}
and  the \emph{stable-spline kernel of order 2}:
\begin{equation}\label{SSkernel0}P^{SS}_{k,j}(\nu) = c\left\{
\begin{array}{cc}
\frac{\mu^{2k}} 2 \left(\mu^j-\frac{\mu^k} 3 \right) & k\geq j\\
\frac{\mu^{2j}} 2 \left(\mu^k-\frac{\mu^j} 3 \right) & k< j
\end{array}\right. ,\ \ \ k,j=1,...,n, \end{equation}
where $\nu = (c,\mu)^T$, $c\geq 0$, $\mu\in[0,1) $. Soon after \cite{SS2010} several other papers appeared where different families of kernels have been introduced \cite{SS2010,SS2011,ChenOL12,SurveyKBsysid}, among which
the diagonal/correlated (DC) kernel
\begin{equation}\label{DCkernel0}P^{DC}_{k,j}(\nu) = c\mu^{(k+j)/2}\rho^{|k-j|},\ \ \ k,j=1,...,n, \end{equation}
\end{subequations}
where $\nu = (c,\mu,\rho)^T$, $c\geq 0$, $\mu\in[0,1) $, $\rho\in (-1,1)$.
As discussed in \cite{ChenetalTAC:14}, and further elaborated upon in \cite{ChiusoCLPCDC2014}, these kernels alone may not well represent impulse responses obtained by linear combination of exponentially decaying functions when the decay rates vary widely; see e.g. Example 2.1 in \cite{ChenetalTAC:14}. For this reason the paper \cite{ChenetalTAC:14} introduces a family of \emph{multiple kernels}, which take the form
\begin{equation}\label{Palpha}
P(\nu) = \sum_{i=1}^m\nu_iP_i,
\end{equation}
where $P_i\in\R^{n\times n}$ are given fixed symmetric and positive semidefinite matrices and the coefficients $\nu_i\geq 0$ ($i=1,\dots,m$) play the role of scale factors.\\
Here, as in \cite{ChenetalTAC:14}, the ``alphabet'' of kernels $P_i$ is chosen from one of the kernels \eqref{nonlinearkernels} over a suitable grid of hyperarameters $(\rho,\mu,c)$. All the kernel choices listed above correspond to an optimization problem \eqref{pr1}--\eqref{fobj} with box-type constraints.\\

\paragraph{Remark}
In our approach, $\sigma^2$ is treated as an optimization variable as suggested in \cite{ChenetalTAC:14}.
As an alternative, the noise variance $\sigma^2$ can be estimated from the data using a high order (and thus low bias) ARX model (linear regressions) as suggested in \cite{Goodwin1992,Ljung:99}; some care needs to be taken to avoid overfitting. In this case only $\nu$, which corresponds to the first $m$ components of $x$, would have to be optimized using the marginal likelihood.

\section{Problem features}\label{sec:problem}

In this section we describe some properties of the optimization problem \eqref{pr1}--\eqref{fobj}. We first need to introduce some notation, defining  the objective function as
\begin{equation}\label{fobjx}
f(x) = f_0(\x) + f_1(\x),
\end{equation}
with
\begin{equation*}
f_0(\x) = Y^T\Sigma(x)^{-1}Y, \ \ \ \ f_1(\x) = \log\det(\Sigma(x)),\ \ \ \forall x\in \R^{m+1}.
\end{equation*}
The $i$--th component of the gradient of $f_0(\x)$ and $f_1(\x)$ can be expressed as
\begin{align}
\label{grad2} \nabla_i f_0(x) &= -Y^T\Sigma(x)^{-1}\frac{\partial \Sigma(x)}{\partial x_i}\Sigma(x)^{-1}Y \\
\label{grad1} \nabla_i f_1(x) &= \mbox{Tr}\left(\Sigma(x)^{-1}\frac{\partial \Sigma(x)}{\partial x_i}\right)
\end{align}
where
\begin{equation*}
\frac{\partial \Sigma(x)}{\partial x_i} = \left\{ \begin{array}{ll} \displaystyle\Phi\frac{\partial P(\nu)}{\partial \nu_i}\Phi^T, & i=1,...,m\\ I_{N-n}, & i = m+1.\end{array}\right.
\end{equation*}
Moreover, the element $(i,j)$ of the Hessian matrix $\nabla^2 f(\x)$, for $i,j=1,...,m+1$, is given by $\nabla^2_{ij}f(\x) = \nabla^2_{ij}f_0(\x)+ \nabla^2_{ij}f_1(\x)$, where
\begin{align*}
\nabla^2_{ij} f_0(\x) \! &= \! Y^T\Sigma(x)^{-1} \! \left( \! \frac{\partial \Sigma(\x)}{\partial x_j}\Sigma(x)^{-1}\frac{\partial \Sigma(\x)}{\partial x_i} \!-\! \frac{\partial^2 \Sigma(\x)}{\partial x_i\partial x_j} \!+\! \frac{\partial \Sigma(\x)}{\partial x_i}\Sigma(x)^{-1}\frac{\partial \Sigma(\x)}{\partial x_j} \! \right) \! \Sigma(x)^{-1} Y, \\
\nabla^2_{ij} f_1(\x) \! &= \! \Tr\left( - \Sigma(x)^{-1}\frac{\partial \Sigma(\x)}{\partial x_j}\Sigma(x)^{-1}\frac{\partial \Sigma(\x)}{\partial x_i} + \Sigma(x)^{-1}\frac{\partial^2 \Sigma(\x)}{\partial x_i\partial x_j}\right).
\end{align*}
When $P(\nu)$ is the multiple kernel \eqref{Palpha}, then $f_0(\x)$ and $f_1(\x)$ are convex and concave, respectively (see \cite{Boyd-Vandeberghe}). In this case, since $\frac{\partial P(\nu)}{ \partial \nu_i}=P_i$, $i=1,\dots,m$, is positive semidefinite, the gradient of the objective function has the following interesting property
\begin{equation}\label{gradpos}
\nabla f_0(\x) \leq 0,\ \ \ \nabla f_1(\x) > 0, \ \ \ \forall x\in \R^m
\end{equation}
when $\sigma^2>0$. The first inequality is straightforward since $\Sigma(x)$ is positive definite and $\Phi P_i\Phi^T $ is positive semidefinite for all $i=1,...,m$, while the second one is a direct consequence of Lemma II.1 in \cite{Lasserre1995}. Moreover, the objective function satisfies
\begin{equation}\label{coercive}
\lim_{t\rightarrow+\infty} f(tx) = +\infty
\end{equation}
for all $x>0$, where $t\in \R$, so that its level sets are bounded, see \cite[\S III.B]{ChenetalTAC:14}.\\
Observe that \eqref{gradpos} and \eqref{coercive} are, in general, not true when the kernel $P(\nu)$ nonlinearly depends on its parameter $\nu$, as in \eqref{nonlinearkernels}; in this case it is not even ensured that $f_0(\x)$ and $f_1(\x)$ are convex and concave, respectively.

\section{Optimization method}\label{sec:opt}

In this section we describe the optimization method we propose to solve \eqref{pr1}. We focus on first order methods based on gradient projection, which are particularly suited when the constraints are simple. The main objection in the use of first order methods is that their convergence rate is, in general, linear. However, introducing some clever choices to define the descent direction, they are able to compute a medium accuracy solution with a small number of iterations.\\
Such acceleration strategies are implemented in the scaled gradient projection (SGP) method \cite{BZZ09}, which applies to any problem of the form
\begin{equation}\label{minfOmega}
\min_{x\in\Omega} f(x),
\end{equation}
where $\Omega\subseteq \R^p$ is a closed convex set, and employs a double scaling of the negative gradient direction through a positive scalar parameter $\alpha_k$ and a positive definite matrix $D_k$, both iteration dependent. The general scheme of SGP is summarized in Algorithm \ref{GPM}. To motivate the introduction of the scaling matrix, one can think, for example, to the Newton's method, which actually scales the gradient direction with the inverse Hessian, while other practical choices for $\alpha_k$ and $D_k$ are described in the following sections.\\
In order to define a descent direction at Step 3, i.e. a vector $\Delta \xk$ such that $\nabla f(\xk)^T \Delta \xk<0$, the projection at Step 2 is computed with respect to the norm induced by the inverse of the scaling matrix $D_k$, i.e. it is defined as
\begin{equation*}
\Pi_{\Omega,D_k^{-1}}(z)= \mbox{arg}\min_{x\in\Omega} (x-z)^T D_k^{-1}(x-z).
\end{equation*}
Thus, even if any positive definite matrix is allowed, the most practical choice for $D_k$ consists in a diagonal matrix with positive diagonal entries.
Once defined the descent direction at Step 3, an Armijo backtracking loop computes the steplength $\lambda_k$ to guarantee the sufficient decrease of the objective function \cite[\S 2.2.1]{Bertsekas}, i.e.
\begin{equation}\label{Armijo}
f(x^{(k)} + \lambda_k\Delta \ve x^{(k)})\leq
           f(x^{(k)})+\beta\lambda_k\nabla f(\ve x^{(k)})^T\Delta \ve x^{(k)}.
\end{equation}

\begin{algorithm}[ht]
\begin{footnotesize}
\caption{Scaled gradient projection (SGP) method}
\label{GPM}
Choose the starting point $\ve x^{(0)}\in \Omega$, set the parameters $\beta, \gamma\in (0,1)$,
$0< \alpha_{min} <\alpha_{max}$, $0<L_{min}<L_{max}$ and fix a positive integer $M$.\\[.2cm]
{\textsc{For}} $k=0,1,2,...$ \textsc{do the following steps:}
\begin{itemize}
\item[]
\begin{AlgorithmSteps}[4]
\item[1] Choose the parameter $\alpha_k \in [\alpha_{min},\alpha_{max}]$ and the diagonal scaling matrix $D_k$ such that $L_{min}\leq(D_k)_{ii}\leq L_{max}$, $i=1,...,p$ ;
\item[2] Projection: $z^{(k)} = \Pi_{\Omega,D_k^{-1}}( x^{(k)}-\alpha_kD_k\nabla f( x^{(k)}))$;
\item[3] Descent direction: $\Delta x^{(k)} = z^{(k)}- x^{(k)}$;
\item[4] Set $\lambda_k = 1$;
\item[5] Backtracking loop:
\begin{ifelse}
\item \textsc{If}
$
f(x^{(k)} + \lambda_k\Delta \ve x^{(k)})\leq
           f(x^{(k)})+\beta\lambda_k\nabla f(\ve x^{(k)})^T\Delta \ve x^{(k)}$
      \textsc{Then} \\  \hspace*{.5cm} go to Step 6;
\item \textsc{Else} \\  \hspace*{.5cm} set $\lambda_k = \gamma \lambda_k$ and go to Step 5.
\item \textsc{Endif}
\end{ifelse}
\item[6] Set $\ve x^{(k+1)} = \ve x^{(k)} + \lambda_k \Delta \ve x^{(k)}$.
\end{AlgorithmSteps}
\end{itemize}
\textsc{End}
\end{footnotesize}
\end{algorithm}
The Armijo condition \eqref{Armijo} is crucial for the proof of the following general convergence result, which can be found in \cite[Theorem 2.1]{BZZ09} (see also \cite[Theorem 4.2]{Bonettini11}).
\begin{Thm}\label{teorema}
Let $\{\ve x^{(k)}\}$ be the sequence generated by applying the SGP algorithm to problem \eqref{minfOmega}. Then, every accumulation point $\ve x^*$ of the sequence $\{\ve x^{(k)}\}$
is a constrained stationary point, that is
$$\nabla f(\ve x^*)^T(\ve x-\ve x^*)\geq 0, \ \ \forall\, \ve x \in \Omega.$$
\end{Thm}
We remark that all the iterates generated by SGP belong to the set $\Omega_0 = \{ x\in \Omega: f(x)\leq f(x^{(0)}) \}$. When $f(x)$ is defined as in \eqref{fobj} and $P(\nu)$ has the form \eqref{Palpha}, we recall that \eqref{coercive} holds: this implies that $\Omega_0$ is bounded and, thus, the sequence $\{\xk\}$ admits at least one limit point.\\
Observe that Theorem \ref{teorema} holds without convexity assumptions and for any bounded choice of the stepsize $\alpha_k$ and scaling matrix $D_k$. This freedom of choice can be exploited to significantly improve the practical performances of SGP. In the following we describe the main strategies for the selection of these parameters.

\subsection{Stepsize selection rules}\label{subsec:BB}

Once a scaling matrix $D_k$ has been defined, a well performing choice of the stepsize parameter is the variant of the Barzilai--Borwein rules proposed in \cite{BZZ09}. The rationale behind this idea consists in computing the stepsize $\alpha_k$ so that the matrix $\alpha_kD_k$ approximates in a quasi--Newton sense the inverse Hessian of the objective function. In practice, $\alpha_k$ is computed as the solution of one of the following minimization problems:
\begin{equation}\label{BB}
 \min_{\alpha\in\mathbb{R}}
\|\alpha D_k r^{(k-1)} - {w}^{(k-1)} \|,\ \ \ \ \
\min_{\alpha \in\mathbb{R}} \|r^{(k-1)} -
(\alpha D_k)^{-1}\ve{w}^{(k-1)} \|,
\end{equation}
where $r^{(k-1)} = {x}^{(k)} - {x}^{(k-1)}$ and
${w}^{(k-1)} = \nabla f({x}^{(k)}) - \nabla f({x}^{(k-1)})$. The solutions of the minimum problems in \eqref{BB} are given by
\begin{equation}
\alpha_k^{{BB1}} = \frac{{{r}^{(k-1)}}^T D_k^{-1} D_k^{-1}
{r}^{(k-1)}}{{{r}^{(k-1)}}^T D_k^{-1} {w}^{(k-1)}}\,,
\quad
\alpha_k^{{BB2}} = \frac{{{r}^{(k-1)}}^T D_k {w}^{(k-1)} }
{{{w}^{(k-1)}}^T D_k D_k {w}^{(k-1)}}\,,
\label{alphaBB}
\end{equation}
and, from the computational point of view, they can be computed in ${\mathcal O}(p)$ operations, where $p$ is the number of variables. Actually, the scalar products $ {{r}^{(k-1)}}^T D_k^{-1} {w}^{(k-1)}$ and ${{r}^{(k-1)}}^T D_k {w}^{(k-1)}$ may be negative, leading to negative values in formula \eqref{BB}. If this occurs, we set $\alpha^{BB1}_k=\alpha_{max}$ and
$\alpha^{BB2}_k=\alpha_{max}$ respectively: this choice is based on the observation that the $(k-1)$-th iterate lies in a region where the objective function might have a negative curvature (if $D_k=I$ and $f$ is convex, both the scalar products are non-negative). Thus, taking a long step along the negative gradient could help to go away from a nonoptimal stationary point.\\
It is well known by the recent literature that the best performances are achieved by adaptively alternating the two rules, with a thresholding to keep the stepsize within the prefixed interval $[\alpha_{min},\alpha_{max}]$ (see Step 1 in Algorithm \ref{GPM}). In our implementation we adopt the alternation strategy detailed below:\\[.1cm]
\hspace*{.3cm}
\begin{footnotesize}
\textsc{if} $ {{r}^{(k-1)}}^T D_k^{-1} {w}^{(k-1)} \leq 0$ \textsc{then}\\
\hspace*{.5cm} $\alpha_k^{(1)} = \alpha_{max}$;\\
\hspace*{.3cm}
\textsc{else}\\
\hspace*{.5cm} $\alpha_k^{(1)} = \min\left\{\alpha_{max}, \max\left\{\alpha_{min},\alpha_k^{BB1}\right\}\right\}$; \\
\hspace*{.3cm}
\textsc{endif}
\\[.1cm]
\hspace*{.3cm}
\textsc{if} $ {{r}^{(k-1)}}^T D_k {w}^{(k-1)}\leq 0$ \textsc{then}\\
\hspace*{.5cm} $\alpha_k^{(2)} = \alpha_{max}$;\\
\hspace*{.3cm}
\textsc{else}\\
\hspace*{.5cm} $\alpha_k^{(2)} = \min\left\{\alpha_{max}, \max\left\{\alpha_{min},\alpha_k^{BB2}\right\}\right\}$; \\
\hspace*{.3cm}
\textsc{endif}
\\[.1cm]
\hspace*{.3cm}
\textsc{if} $ {\alpha_k^{(2)}}/{\alpha_k^{(1)}} \le
\tau_k$ \textsc{then}\\
\hspace*{.5cm} $\alpha_k = \min\left\{\alpha_j^{(2)}, \
j=\max\left\{1,k-M_{\alpha}\right\},\dots,k\right\}$; \ $\tau_{k+1} = \tau_{k}\cdot 0.9$;\\
\hspace*{.3cm}
\textsc{else}\\
\hspace*{.5cm} $\alpha_k = \alpha_k^{(1)}$; \ $\tau_{k+1} = \tau_{k}\cdot 1.1$;\\
\hspace*{.3cm}
\textsc{endif}
\end{footnotesize}
\\[.1cm]
where $M_{\alpha}$ is a prefixed non-negative integer and
$\tau_1\in (0,1)$. The alternating rule described above has been proposed for unconstrained, strictly convex quadratic problems in \cite{Frassoldati-etal-2008}, where the authors investigate the related theoretical properties and numerically show that this alternation of the two BB rules allows to better capture the spectral properties of the Hessian matrix. Successively, an adaptation of the alternating rule in \cite{Frassoldati-etal-2008} has been proposed in \cite{BZZ09} and employed also in several applications of SGP to different convex, nonlinear, constrained problems \cite{Bonettini11,BonettiniCornelioPrato2013,Bonettinietal2013,Bonettini2010}. In this paper we adopt the same rule also for the nonlinear, nonconvex, constrained problem described in Section \ref{sec:problem}.

\subsection{Choice of the scaling matrix}\label{sec:scaling}

Unlike the stepsize selection rules, the scaling matrix choice is strictly related to the specific structure of problem \eqref{minfOmega} and it depends on both the objective function and the constraints. In particular, the constraints of problem \eqref{pr1} are lower bounds when $P(\nu)$ is the multiple kernel \eqref{Palpha} or box constraints when the kernels \eqref{nonlinearkernels} are selected.\\
In this section we review the split gradient idea described in \cite{Bertero-etal-08,Lanteri-etal-02} for lower bound constraints and we extend such approach to general box constraints. To introduce the split gradient idea, we consider first the non-negatively constrained problem
\begin{equation}\label{minfnonneg}
\min_{\x\geq 0} f(\x)
\end{equation}
whose first order optimality conditions are given by
\begin{equation}\label{equality}
x\nabla f(x) = 0; \quad \ x\geq 0; \quad \nabla f(x) \geq 0,
\end{equation}
where the equality and inequalities are componentwise. If the gradient of $f(\x)$ admits a decomposition like the following one
\begin{equation}\label{SplitGrad}\nabla f(x) = V(x)-U(x)\ \mbox{ with } V(x) > 0, \ U(x)\geq 0, \end{equation}
then equality \eqref{equality} writes also as the fixed point equation $x = xU(x)/V(x)$. This formulation is related to the corresponding fixed point method
\begin{equation}\label{mult}
x^{(k+1)} = x^{(k)}\frac{U(x^{(k)})}{V(x^{(k)})},
\end{equation}
whose convergence properties are not well studied, but which has the capability to preserve positivity when the initial point is positive and $U(x)>0$ whenever $x>0$. Several methods in signal and image processing (e.g. Lucy-Richardson/expectation minimization \cite{Shepp1982}, iterative space reconstruction algorithm \cite{Daube-Witherspoon1986}) and statistical learning (Lee-Seung algorithm for non-negative matrix factorization \cite{LeeSeung1999}) actually have exactly this multiplicative form (see also \cite{Hageretal09,Merritt2005}).\\
With a simple algebra the multiplicative method \eqref{mult} results in
\begin{equation*}
x^{(k+1)} = \ x^{(k)}\frac{U(x^{(k)})-V(x^{(k)})+ V(x^{(k)})}{V(x^{(k)})}=\ x^{(k)} - \frac{x^{(k)}}{V(x^{(k)})}\nabla f(x^{(k)}),
\end{equation*}
which corresponds to a scaled gradient iteration. These considerations suggest to define the scaling matrix for SGP as
\begin{equation}\label{scaling0}(D_k)_{ii} = \min\left( \max\left(L_{min}, \frac{x_i^{(k)}}{V_i(x^{(k)})}\right),L_{max}\right) . \end{equation}
More in general, for lower bound constraints $\x\geq l$, $l\in \R^p$, the following scaling matrix
\begin{equation}\label{scaling}(D_k)_{ii} = \min\left( \max\left(L_{min}, \frac{x_i^{(k)}-l_i}{V_i(x^{(k)})}\right),L_{max}\right)  \end{equation}
can be motivated using similar arguments as above.\\
This choice of the scaling matrix, combined with a suitable choice of the stepsize $\alpha_k$, leads the SGP method to very good performances on ill-posed/ill-conditioned inverse problems approached by the Bayesian paradigm as convex, non-negatively constrained optimization problems \cite{BZZ09,Prato2012,Zanella2009,Zanella2013b}.\\
We propose to use the scaling \eqref{scaling} also on problem \eqref{pr1} with the multiple kernel \eqref{Palpha}, even if the objective function is nonconvex. In this case, recalling \eqref{gradpos}, the gradient of the objective function has the natural decomposition \eqref{SplitGrad} with
\begin{equation}\label{graddec1}
V(\x) = \nabla f_1(\x) \ \ \ \mbox{ and } \ \ \ U(\x) = -\nabla f_0(\x).
\end{equation}

\subsubsection{Gradient splitting strategy for box constraints}

In this section we consider a box constrained problem
\begin{equation}\label{minfbox}
\min_{l\leq \x\leq u} f(\x)
\end{equation}
where $l,u\in \R^p\cup\{+\infty,-\infty\}$ ($l_i=-\infty$, $u_i=+\infty$ means that $x_i$ is unbounded below or above respectively), and we propose a scaling strategy also for this case. Driven by the considerations made in the previous section, the generalization to box constraints consists in finding a positive diagonal scaling matrix $D_k$ such that $\xk - D_k\nabla f(\xk)$ is feasible, i.e.
\begin{equation}\nonumber
l_i\leq \xk_i - (D_k)_{ii}\nabla_i f(\xk)\leq u_i, \ \ \ i=1,\dots,p.
\end{equation}
Then, to design an appropriate scaling, we should consider the sign of the gradient at the current iterate to devise which constraints could be violated taking a step along the negative gradient direction. To this end, we define the following sets of indices
\begin{eqnarray*}
{\mathcal I}_1 = \{i: l_i > -\infty\mbox{ and } u_i < +\infty\} &, &\ \ {\mathcal I}_2 = \{i: l_i = -\infty\mbox{ and } u_i < +\infty\},\\
{\mathcal I}_3 = \{i: l_i > -\infty\mbox{ and } u_i = +\infty\} &, &\ \ {\mathcal I}_4 = \{i: l_i = -\infty\mbox{ and } u_i = +\infty\},
\end{eqnarray*}
to identify which variables are bounded below and/or above and which are unbounded. Then, we define the following vector
 \begin{equation}\label{scaling-bounds1}
\tilde d_i(\xk) = \left\{ \begin{array}{ll}
\displaystyle \frac{u_i-\xk_i}{U_i(\xk)} & \mbox{ if } i\in{\mathcal I}_1\mbox{ and }\nabla_i f(\xk) \leq 0   \mbox{ or } i\in {\mathcal I}_2\vspace{0.1cm}\\
\displaystyle \frac{\xk_i-l_i}{V_i(\xk)} & \mbox{ if } i\in{\mathcal I}_1\mbox{ and }\nabla_i f(\xk) > 0 \mbox{ or } i\in {\mathcal I}_3\vspace{0.1cm}\\
\displaystyle 1                          & \mbox{ if } i\in{\mathcal I}_4
\end{array}
\right.
\end{equation}
based on a gradient decomposition of the form
\begin{equation}\label{graddec2}
\nabla f(\x) = V(\x)-U(\x), \ \ V(x) > 0,\ U_i(x)>0.
\end{equation}
Indeed, $\nabla_i f(\xk) \leq 0$ implies $0<V_i(\xk)\leq U_i(\xk)$ and, as a consequence, $\xk_i \leq \xk_i -\tilde d_i (\xk) \nabla_if(\xk) \leq u_i$. On the other side, $\nabla_i f(\xk) > 0$ if and only if $V_i(\xk)\geq U_i(\xk)> 0$, which yields $ l_i\leq \xk_i -\tilde d_i(\xk) \nabla_if(\xk) \leq \xk_i$.\\ Finally, the diagonal entries of the scaling matrix are defined as
\begin{equation}\label{scaling-bounds2}(D_k)_{ii} = \min\left( \max\left(L_{min}, \tilde{d}_i(\xk)\right),L_{max}\right) . \end{equation}
For an objective function of the form $f(x) = f_0(x) + f_1(x)$ a possible general rule to define $U(\xk)$ and $V(\xk)$ in \eqref{scaling-bounds1} can be devised in the following way.
When $\nabla_i f(\x)> 0$, then $\nabla_i f_1(\x) > -\nabla_i f_0(\x)$ and we define
\begin{equation}\label{Vbounds}
\begin{array}{l}
V_i(\x) = \left\{ \begin{array}{lllll}
\nabla_i f_0(\x) & \mbox{ if }  \nabla_i f_1(\x) < 0 \\
\nabla_i f_1(\x) & \mbox{ if }   \nabla_i f_1(\x) \geq 0 \mbox{ and } \nabla_i f_0(\x) < 0\\
\nabla_i f(\x) + \zeta   & \mbox{ otherwise}
\end{array}
\right.,\vspace{0.2cm}\\ U_i(\x) = V_i(\x)-\nabla_i f(\x)\end{array}
\end{equation}
for some $\zeta > 0$.
Similarly, when $\nabla_i f(\x) \leq 0$, then $\nabla_i f_1(\x) \leq -\nabla_i f_0(\x)$ and we set
\begin{equation}\label{Ubounds}
\begin{array}{l}
U_i(\x) = \left\{ \begin{array}{lllll}
-\nabla_i f_1(\x)   & \mbox{ if }  \nabla_i f_0(\x) > 0 \\
-\nabla_i f_0(\x)   & \mbox{ if }  \nabla_i f_0(\x) < 0 \mbox{ and } \nabla_i f_1(x) > 0\\
\zeta -\nabla_if(\x)&            \mbox{ otherwise}
\end{array}
\right.,\vspace{0.2cm} \\ V_i(\x) = \nabla_i f(\x) + U_i(\x).\end{array}
\end{equation}
It is easy to verify that definitions \eqref{Vbounds} and \eqref{Ubounds} lead to a gradient decomposition with the property \eqref{graddec2}. Moreover, this choice of the scaling matrix reduces to \eqref{scaling} in presence of lower bounds only.\\
We adopt the scaling strategy \eqref{scaling-bounds2} associated to the decomposition \eqref{Vbounds}--\eqref{Ubounds} for problem \eqref{pr1} when the kernel is given by \eqref{nonlinearkernels}.

\subsection{Algorithm implementation and complexity}\label{sec:implementation}

Each SGP iteration requires the objective function \eqref{fobjx} and gradient \eqref{grad2}--\eqref{grad1} at the current point $\x^{(k)} = ({\nu^{(k)}}^T,\sigma_k^2)^T$, which is the more relevant computational burden of the whole algorithm. If the Armijo condition \eqref{Armijo} is not satisfied with $\lambda_k=1$, more function evaluations are needed.\\Thus, the practical performances of the algorithm also relies on the implementation of the gradient and objective function computation. On the other side we should take into account the severe ill-conditioning possibly affecting the matrices $P(\nu^{(k)})$ and $\Sigma(\xk)$. In our implementation we implicitly assume that $n\ll N$, which is quite realistic, and we devise an algorithm for the computation of $f(\xk)$ and $\nabla f(\xk)$ with complexity ${\mathcal O}(n^3)$ which is detailed below.\\
We consider the approach proposed in \cite{ChenLjung13} for objective function and gradient evaluations, which is based on the Cholesky factorization of $P(\nu^{(k)})={\mathcal L}_k{\mathcal L}_k^T$, at a cost of ${\mathcal O}(n^3)$. Then, the Cholesky factorization of the matrix $\sigma^2_k I_n+ {\mathcal L}_k^T\Phi^T \Phi{\mathcal L}_k = S_kS_k^T$ is also computed. Finally, the objective function is evaluated with the formula
\begin{equation}\label{fobj2}
f(\xk) = (\|Y\|^2 -\|S_k^{-1}{\mathcal L}_k^T\Phi^T  Y \|^2)/\sigma_k^2 + (N-n)\log \sigma_k^2 + 2\log |S_k|.
\end{equation}
The Cholesky factors $S_k$ and ${\mathcal L}_k$ can be reused for the computation of $\Sigma(\xk)^{-1}$ and, then, of the gradient as follows. Omitting for simplicity the dependency of $P(\nu^{(k)})$ and $\Sigma(\xk)$ from $\nu^{(k)}$ and $\xk$
and applying the Sherman-Morrison-Woodbury formula we obtain
\begin{equation*}
\Sigma^{-1} = (\sigma_k^2I_{N-n} + \Phi P\Phi^T )^{-1} = \frac 1 {\sigma_k^2}I_{N-n}-\frac 1 {\sigma_k^2}\Phi(\sigma_k^2P^{-1}+\Phi^T \Phi)^{-1}\Phi^T .
\end{equation*}
Finally, by observing that
\begin{equation*}
 (\sigma_k^2 P^{-1} +\Phi^T \Phi)^{-1} = (\sigma_k^2{\mathcal L}_k^{-T}{\mathcal L}_k^{-1} + \Phi^T \Phi)^{-1} = {\mathcal L}_kS_k^{-T}S_k^{-1}{\mathcal L}^T_k
\end{equation*}
it follows that
\begin{equation}\label{Sigmacompute}
\Sigma^{-1} = \frac 1 {\sigma_k^2}I_{N-n}-\frac 1 {\sigma_k^2}\Phi{\mathcal L}_kS_k^{-T}S_k^{-1}{\mathcal L}_k^T\Phi^T.
\end{equation}
Taking into account of \eqref{Sigmacompute}, if we set
\begin{equation}\label{setting}
\tilde \Phi = \Phi^T\Phi,\ \ \tilde Y = \Phi^T Y,\ \  Z_k = {\mathcal L}_kS_k^{-T}S_k^{-1}{\mathcal L}_k^T, \ \ M_k = \Phi^T\Sigma^{-1}\Phi = \tilde\Phi-\tilde\Phi Z_k \tilde\Phi,
\end{equation}
then \eqref{grad2} can be computed as
\begin{subequations}\label{Ucompute}
\begin{equation}\label{Ucomputea}
\nabla_if_0(\xk) =  q^T \frac{\partial P}{\partial \nu_i}q, \ \ \mbox{ with } q =\Phi\Sigma^{-1}Y = \frac{1}{\sigma_k^2}(I_n-\tilde\Phi Z_k)\tilde Y
\end{equation}
for $i=1,...,m$ and
\begin{eqnarray}
\nabla_{m+1}f_0(\xk) &=&  \|\Sigma^{-1}Y\|^2 \nonumber\\
&=& \|Y\|^2/\sigma_k^4 -2 Y^T\Phi{\mathcal L}_kS^{-T}_kS_k^{-1}{\mathcal L}_k^T\Phi^TY/\sigma_k^4 + \nonumber\\
& &+Y^T\Phi{\mathcal L}_kS_k^{-T}S_k^{-1}({\mathcal L}_k^T \Phi^T\Phi{\mathcal L}_k)S^{-T}_kS_k^{-1}{\mathcal L}_k^T\Phi^TY/\sigma_k^4\nonumber\\
&=& \|Y\|^2/\sigma_k^4 -2 Y^T\Phi{\mathcal L}_kS^{-T}_kS_k^{-1}{\mathcal L}_k^T\Phi^TY/\sigma_k^4 +\nonumber\\
& &+Y^T\Phi{\mathcal L}_kS_k^{-T}S_k^{-1}(S_kS_k^T-\sigma_k^2I_n)S^{-T}_kS_k^{-1}{\mathcal L}_k^T\Phi^TY/\sigma_k^4\nonumber\\
&=&\|Y\|^2/\sigma_k^4 - \|S_k^{-1}{\mathcal L}_k^T\tilde Y\|^2/\sigma_k^4 - \|S_k^{-T}S_k^{-1}{\mathcal L}_k^T\tilde Y\|^2/\sigma_k^2.\label{Ucomputeb}
\end{eqnarray}
\end{subequations}
On the other side, recalling \eqref{grad1}, we have
\begin{subequations}\label{Vcompute}
\begin{eqnarray}
\nabla_if_1(\xk) &= &\mbox{Tr}\left( \Sigma^{-1} \Phi\frac{\partial P}{\partial \nu_i}\Phi^T \right)
= \mbox{Tr}\left( \Phi^T \Sigma^{-1} \Phi\frac{\partial P}{\partial \nu_i}\right)\nonumber\\
&= &\frac 1 {\sigma_k^2} \mbox{Tr}\left( M_k \frac{\partial P}{\partial \nu_i}\right), \ \ \ i=1,...,m \label{Vcomputea}
\end{eqnarray}
\begin{eqnarray}
\nabla_{m+1}f_1(\xk) &=& \Tr(I_{N-n} - \Phi{\mathcal L}_kS_k^{-T}S_k^{-1}{\mathcal L}_k\Phi^T)/{\sigma_k^2} \nonumber\\
&=&  \left( \Tr(I_{N-n}) + \Tr({\mathcal L}_k\Phi^T \Phi{\mathcal L}_kS_k^{-T}S_k^{-1})\right)/ \sigma_k^2\nonumber\\
&=&  \left( \Tr(I_{N-n}) + \Tr((S_kS_k^T-\sigma_k^2I_n)S_k^{-T}S_k^{-1})\right) / \sigma_k^2\nonumber\\ &=& (N-2n)/\sigma_k^2+\Tr(S_k^{-T}S_k^{-1}).\label{Vcomputeb}
\end{eqnarray}
\end{subequations}
The main difference between our approach for gradient computation and the analogous one described in \cite[Section 5]{ChenLjung13} is that we avoid to explicitly compute the matrix $P^{-1} ={\mathcal L}_k^{-T}{\mathcal L}_k^{-1} $, which is very ill-conditioned.\\
The previous formulae for gradient computation clearly hold when $P(\nu^{(k)})$ does not reduce to the zero matrix: since the latter case can occur at some iteration $k$, for sake of completeness we report the whole procedure in Algorithm \ref{algo:fobj_gradobj}.

\begin{algorithm}[ht]
\begin{footnotesize}
\caption{Objective function and gradient evaluation}
\label{algo:fobj_gradobj}
Preprocessing: compute $\|Y\|^2$, $\tilde\Phi = \Phi^T \Phi$ and $\tilde Y = \Phi^T  Y$ .\\[.2cm]
{\textsc{For any}} $\xk=({\nu^{(k)}}^T,\sigma_k^2)^T$, $k=1,2,...$  \textsc{do the following steps:}
\begin{itemize}
\item[]
\begin{AlgorithmSteps}[4]
\item[1] Compute $P(\nu^{(k)})$.
\item[2] \begin{ifelse}
\item \hspace*{-.8cm} \textsc{If} $P(\nu^{(k)})\neq 0$ \textsc{Then}
\hspace*{1.5cm} \item[2.1] Compute the Cholesky factorization $P(\nu^{(k)})={\mathcal L}_k{\mathcal L}_k^T$;
\hspace*{1.5cm} \item[2.2] Compute $Q_k =\sigma_k^2I_n+{\mathcal L}_k^T\tilde\Phi {\mathcal L}_k$;
\hspace*{1.5cm} \item[2.3] Compute the Cholesky factorization $Q_k=S_kS_k^T$, $Z_k$ and $M_k$ as in \eqref{setting};
\hspace*{1.5cm} \item[2.4] Compute $f(\xk)$ by formula \eqref{fobj2};
\hspace*{1.5cm} \item[2.5] Compute $\nabla_i f_0(\xk)$ and $\nabla_i f_1(\xk)$ by means of \eqref{Ucompute} and \eqref{Vcompute} for $i=1,...,m+1$.
\item \hspace*{-.8cm} \textsc{Else}
\hspace*{1.5cm} \item[2.6] Compute $f(\xk) = \frac{\|Y\|^2}{\sigma_k^2} + (N-n)\log(\sigma_k^2)$;
\hspace*{1.5cm} \item[2.7] Compute $\nabla_if_1(\xk) = \frac{1}{\sigma_k^2}\mbox{Tr}\left( \tilde \Phi\frac{\partial P(\nu^{(k)})}{\partial x_i} \right)$ and $\nabla_if_0(\xk) = -\frac{1}{\sigma_k^4} \tilde Y^T\frac{\partial P(\nu^{(k)})}{\partial x_i}\tilde  Y$ for $i=1,...,m$; $\nabla_{m+1} f_1(\xk) = (N-n)/\sigma_k^2$; $\nabla_{m+1}f_0(\xk) = -\|Y\|^2/\sigma_k^4$;
\item \hspace*{-.8cm} \textsc{Endif}
\end{ifelse}
\item[3] Compute $\nabla f(\xk) = \nabla f_0(\xk)+\nabla f_1(\xk)$.
\end{AlgorithmSteps}
\end{itemize}
\textsc{End}
\end{footnotesize}
\end{algorithm}

\paragraph{Remark} The computation of the Hessian matrix can also be performed with a complexity of $\mathcal O (n^3)$, without need of further factorizations but with at least $m$ additional matrix-matrix products of size $n\times n$, as detailed in the following.
Developing the formulae for the entries of the Hessian matrix given in Section \ref{sec:problem}, for $i,j=1,...,m$, we can set $\nabla^2_{ij} f_0(x) = a_{ij}-b_{ij}+a_{ji}$, where
\begin{eqnarray*}
a_{ij} &=& Y^T\Sigma^{-1}\Phi\frac{\partial P}{\partial \nu_j} \Phi^T\Sigma^{-1}\Phi\frac{\partial P}{\partial \nu_i}\Phi^T\Sigma^{-1}Y= q^T\frac{\partial P}{\partial \nu_j} M_k\frac{\partial P}{\partial \nu_i}q\\
b_{ij} &=& Y^T\Sigma^{-1} \frac{\partial^2 \Sigma}{\partial x_i\partial x_j} \Sigma^{-1} Y = q^T \frac{\partial^2 P}{\partial \nu_i\partial \nu_j}q
\end{eqnarray*}
with $q$ defined as in \eqref{Ucomputea}. Moreover we have
\begin{eqnarray*}
\nabla^2_{i,m+1}f_0(x^{(k)})&=& \tilde q\frac{\partial P}{\partial \nu_i} q,\ \ \ i=1,...,m\\
\nabla_{m+1,m+1}^2f_1(x^{(k)}) &=& Y^T\Sigma^{-3} Y \\ &=& \frac{1}{\sigma_k^6}\left(\|Y\|^2 -3\tilde Y Z_k \tilde Y + 3 \tilde Y^TZ_k \tilde \Phi Z_k\tilde Y - \tilde Y^T Z_k \tilde \Phi Z_k \tilde \Phi Z_k \tilde Y\right)
\end{eqnarray*}
where $\tilde q = \Phi^T\Sigma^{-2}Y = (I_n-\tilde \Phi Z_k)^2 \tilde Y/\sigma_k^4$.
As concerns the Hessian of $f_1$, exploiting the matrix trace properties,  for $i,j=1,...,m$ we have $\nabla_{ij}^2f_1(x) = \bij_{ij}-\aij_{ij}$, where
\begin{eqnarray*}
\aij_{ij} &=& \Tr\left( \Sigma^{-1}\frac{\partial \Sigma}{\partial x_j}\Sigma^{-1}\frac{\partial \Sigma}{\partial x_i} \right)
= \Tr\left(M_k\frac{\partial P}{\partial \nu_j}M_k\frac{\partial P}{\partial \nu_i}\right)\\
\bij_{ij} &=& \Tr\left(\Sigma^{-1}\frac{\partial^2 \Sigma}{\partial x_i\partial x_j}\right)= \Tr\left( M_k \frac{\partial^2 P}{\partial \nu_i\partial \nu_j}\right)
\end{eqnarray*}
with $M_k$ defined as in \eqref{setting}, and
\begin{eqnarray*}
\nabla_{i,m+1}^2f_1(x^{(k)}) &=& \Tr\left( \Phi^T \Sigma^{-2} \Phi \frac{\partial P}{\partial \nu_i}\right)= \Tr\left( \tilde M_k\frac{\partial P}{\partial \nu_i}\right),\ \ \ i=1,...,m\\
\nabla_{m+1,m+1}^2f_1(x^{(k)}) &=& \Tr(\Sigma^{-2}) = \frac 1{\sigma_k^4}\left( N-n - 2 \Tr(\tilde \Phi Z_k) + \Tr(\tilde\Phi Z_k\tilde\Phi Z_k) \right)
\end{eqnarray*}
with $\tilde\Phi$, $Z_k$ defined as in \eqref{setting} and $\tilde M_k = (I_n-\tilde\Phi Z_k)\tilde\Phi/\sigma_k^4$.
Observe that $\aij_{ij}$ requires the explicit computation of the matrices $M_k\frac{\partial P}{\partial \nu_i}$, $i=1,...,m$, with a complexity of ${\mathcal O}(mn^3) $. For the multiple kernel \eqref{Palpha}, where $m$ typically is of order of tenths, the Hessian computation is a quite expensive task. It is worth stressing that the computation of the matrix product $M_k\frac{\partial P}{\partial \nu_i}$ is not needed for \eqref{Vcomputea}, since the well known formula $\Tr(AB)=\mbox{vec}(A)^T\mbox{vec}(B)$, where $\mbox{vec}(\cdot)$ indicates the vectorization of a matrix by stacking its elements columnwise, can be applied.
\section{Numerical experience}\label{sec:num}

We consider the test sets described in \cite[Section V.A]{ChenetalTAC:14}, containing 1000 simulated data records $\{y(t),u(t)\}_{t=1}^N$:
\begin{itemize}
\item D1: $N=210$, output SNR = 10;
\item D2: $N=210$, output SNR = 1;
\item D3: $N=500$, output SNR = 10;
\item D4: $N=500$, output SNR = 1.
\end{itemize}
The estimated model order is set to $n=100$ for all simulations.\\
We consider two sets of test problems. In the first one, we choose $P(x)$ as the multiple kernel \eqref{Palpha}, where the `basis' matrices $P_i$ are chosen as follows:
\begin{itemize}
\item[$\bullet$] [DC-M]: $P_i = P^{DC}(1,\mu_i,\rho_i)$ where $P^{DC}$ is the DC kernel defined in \eqref{DCkernel0}
and $(\mu_i,\rho_i)$ are points of the grid
\begin{equation*}\{0.1i:1\leq i\leq 9\}\times\{-0.95,-0.65,-0.35,0.35,0.65,0.95\}\end{equation*}
so that $m=54$;
\item[$\bullet$] [TCSS-M]: $P_i = P^{TC}(1,\mu_i^{TC})$, $i=1,...,21$ where $P^{TC}$ is defined in \eqref{TCkernel0} and $\mu_i^{TC}\in\{0.05i: 2\leq i\leq 15\}\cup\{0.81+0.02i:0\leq i\leq 6\}$, $P_{21+i} = P^{SS}(1,\mu_i^{SS})$, $i=1,...,8$ where $P^{SS}$ is defined in \eqref{SSkernel0} and $\mu_i^{SS} \in \{0.8+0.02i:0\leq i\leq 7\}$. In this case we have $m=29$.
\end{itemize}
The matrices $P_i$, $i=1,...,m$ are extremely ill-conditioned: indeed, the average condition number is about $ 10^n$.
As concerns the choice of $m$, we performed several tests also with finer grids and we observed similar behaviours of the algorithms with no significant improvements in the quality of the estimated impulse response coefficients.\\
In the second set of problems, we consider the following cases:
\begin{itemize}
\item[$\bullet$][DC] $P(x)$ is the DC kernel \eqref{DCkernel0} with $x = (c,\mu,\rho)^T$, where $c\geq 0$, $0.72\leq \mu\leq 0.99$, $-0.99\leq\rho\leq 0.99$;
\item[$\bullet$][TC] $P(x)$ is the TC kernel \eqref{TCkernel0} with $x = (c,\mu)^T$, where $c\geq 0$, $0.7\leq \mu\leq 0.99$;
\item[$\bullet$][SS] $P(x)$ is the SS kernel \eqref{SSkernel0} with $x = (c,\mu)^T$, where $c\geq 0$, $0.7\leq \mu\leq 0.99$.
\end{itemize}
We choose the lower bounds on the `$\mu$' variable according to \cite{ChenLjung13}, with the aim to impose a reasonable upper bound to the condition number of $P(x)$.\\
The quality of the estimated models $\hat \theta$ is evaluated by the coefficient
\begin{equation}\label{fitdef}
 W(\hat \theta) = 100\left(1-\sqrt{\frac{\sum_{i=1}^n |\theta_i^*-\hat\theta_i|^2}{\sum_{i=1}^n |\theta_i^*-\bar\theta|^2}}\right),\ \ \
\bar\theta = \frac 1 n\sum_{i=1}^n\theta_i^*,
\end{equation}
where $\theta^*_i$ are the true impulse response coefficients and $\hat \theta_i$ the estimated ones computed by formula \eqref{thetahat}.

\subsection{SGP parameters setting}\label{subsec:SGPparams}

The SGP parameters have been set as follows: $\beta = 10^{-4}$, $\gamma = 0.4$, $\alpha_{min} = 10^{-7}$, $\alpha_{max} = 10^2$, $L_{min}=\zeta=10^{-5}$, $L_{max}=10^{10}$, $M_{\alpha}=3$, $\tau = 0.5$. The initial point $x^{(0)}$ is the vector of all ones for the multiple kernels, while we set $x^{(0)} = (0.5,0.5,0.8,0.5)^T$ for the DC kernel and $x^{(0)} = (0.5,0.8,0.5)^T$ for the TC and $SS$ kernels. The initial stepsize $\alpha_0$ is set to 1.\\
Since SGP is a projection method, it can occur that some of the iterates lay on the boundary of the feasible set. This may create some trouble, since for $x_{m+1}=\sigma^2=0$ the matrix $\Sigma(x)$ in \eqref{Sigma} may become singular. For these reasons, we constrain the $(m+1)$-variable to be greater or equal to some positive constant.
Then, we actually consider a problem of the form \eqref{minfbox}
where $l\in\R^{m+1}$, $u\in \R^{m+1}\cup \{+\infty\}$, with $l_{m+1} = 10^{-2}$, $u_{m+1}= + \infty$, and $l_i = 0$, $u_i = +\infty$, $i=1,...,m$ for the multiple kernels DC-M and TCSS, $l_1 = 0$, $l_2 = 0.72$, $l_3 = -0.99$, $u_1=+\infty$, $u_2=0.99$, $u_3 = 0.99$ for the kernel DC and $l_1 = 0$, $l_2 = 0.7$, $u_1=+\infty$, $u_2=0.99$ for the kernels TC and SS.\\
We experimentally observed that the constraint on $\sigma^2$ is never active at the solution of \eqref{minfbox} with $l_{m+1}=10^{-2}$: we experienced also smaller values, down to $10^{-8}$, but we did not observe significant differences in the results and in the algorithms performance. As an alternative, this lower bound can be safely set to a fraction (say between one  tenth to one hundredth) of a preliminary estimate of the noise variance which can be obtained, for instance, as discussed in \cite{SS2010,SS2011}.\\
We include also the following stopping criterion for the iterates:
\begin{equation}\label{stopcrit}
f(\xk)-f(\xkk) < \tau |f(\xkk)|
\end{equation}
with $\tau= 10^{-9}$. Indeed, we experienced different values of $\tau$, ranging from $10^{-11}$ to $10^{-7}$ and we observed that no significant improvements in accuracy are obtained with smaller tolerance values. A maximum number of 5000 iterations is also imposed.

\subsection{Scaling matrix impact}

\begin{figure}[ht]
\begin{center}
\begin{tabular}{ccc}
  \includegraphics[width=0.3\linewidth]{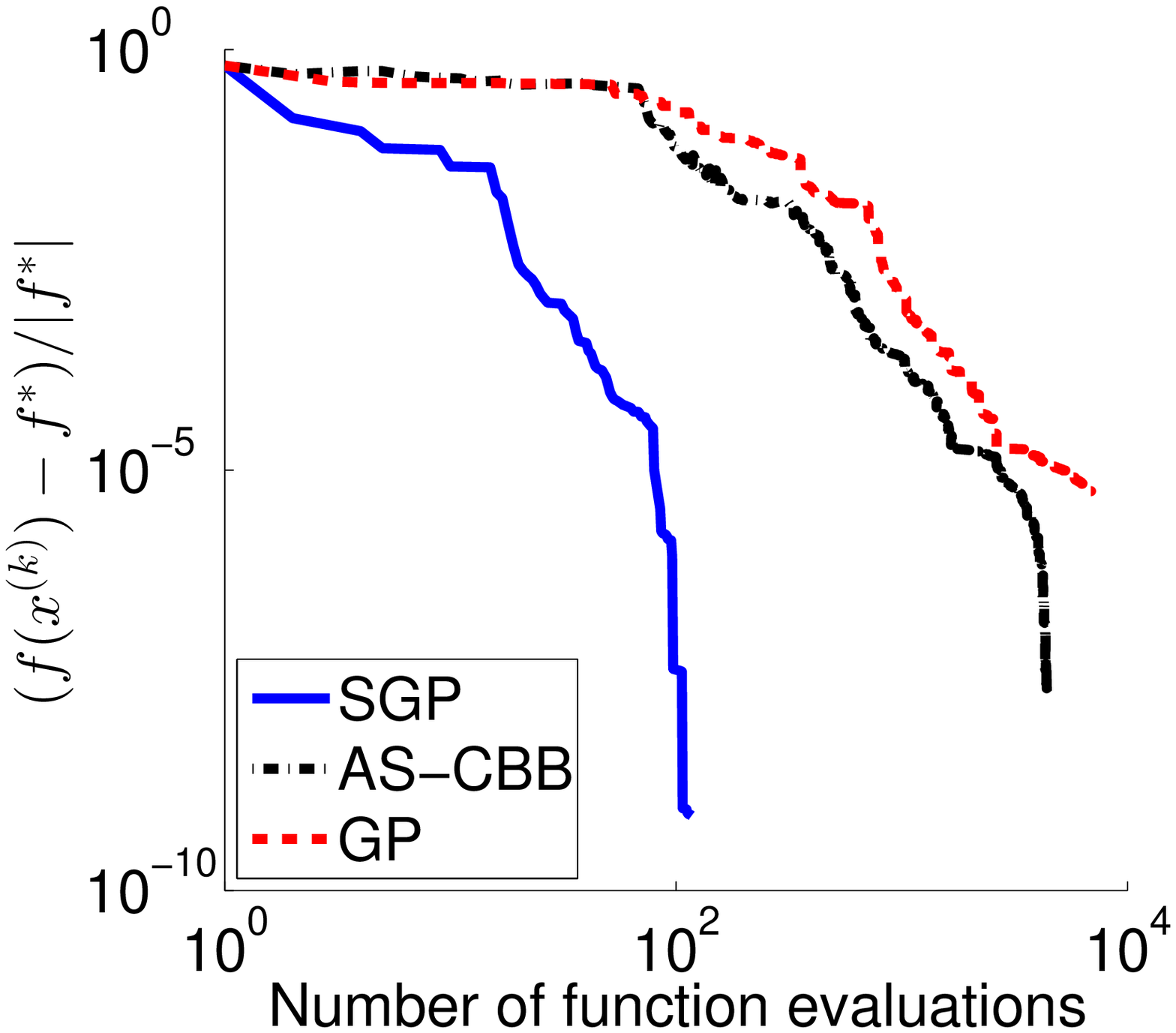} &
  \includegraphics[width=0.3\linewidth]{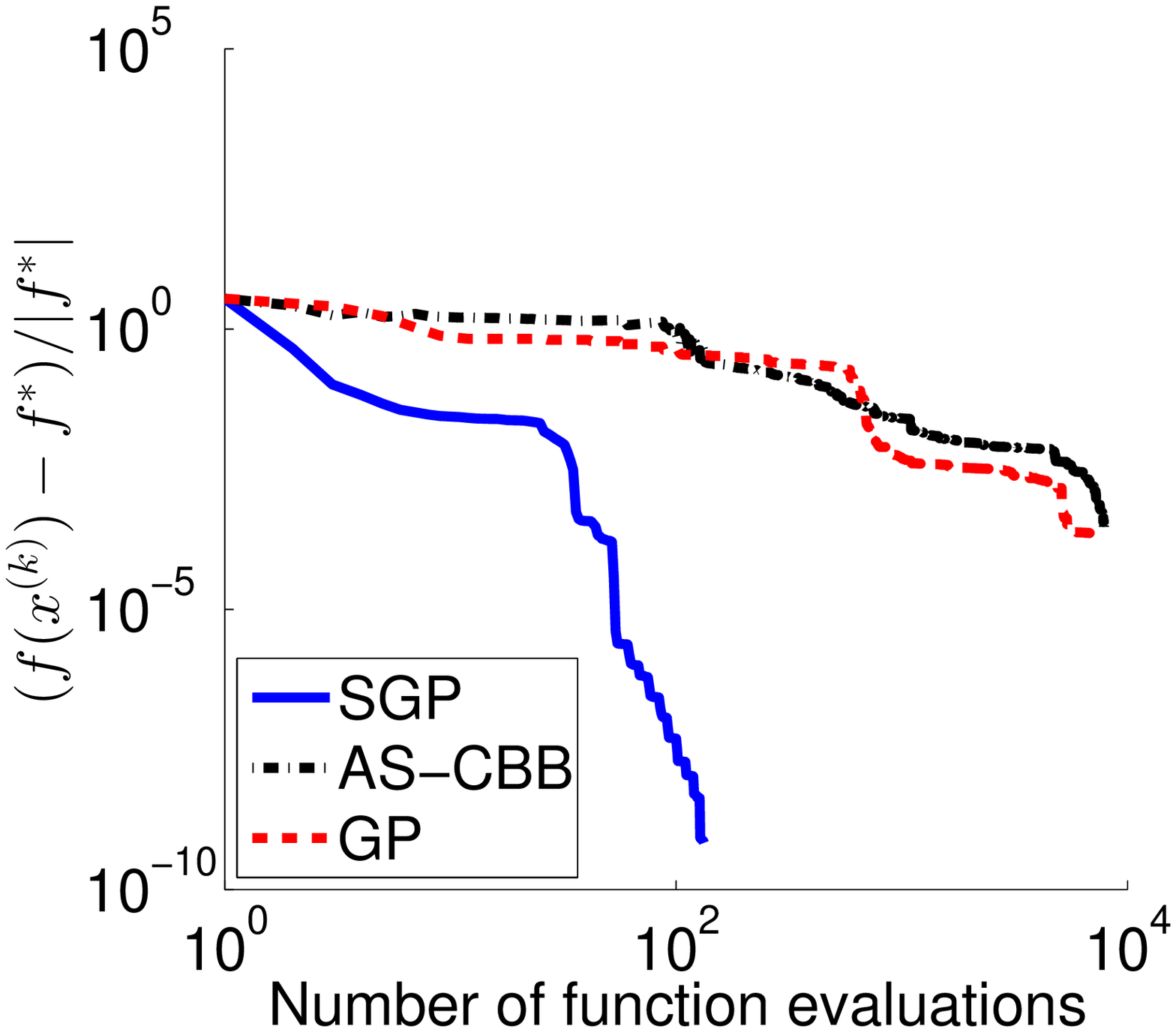} &
	\includegraphics[width=0.3\linewidth]{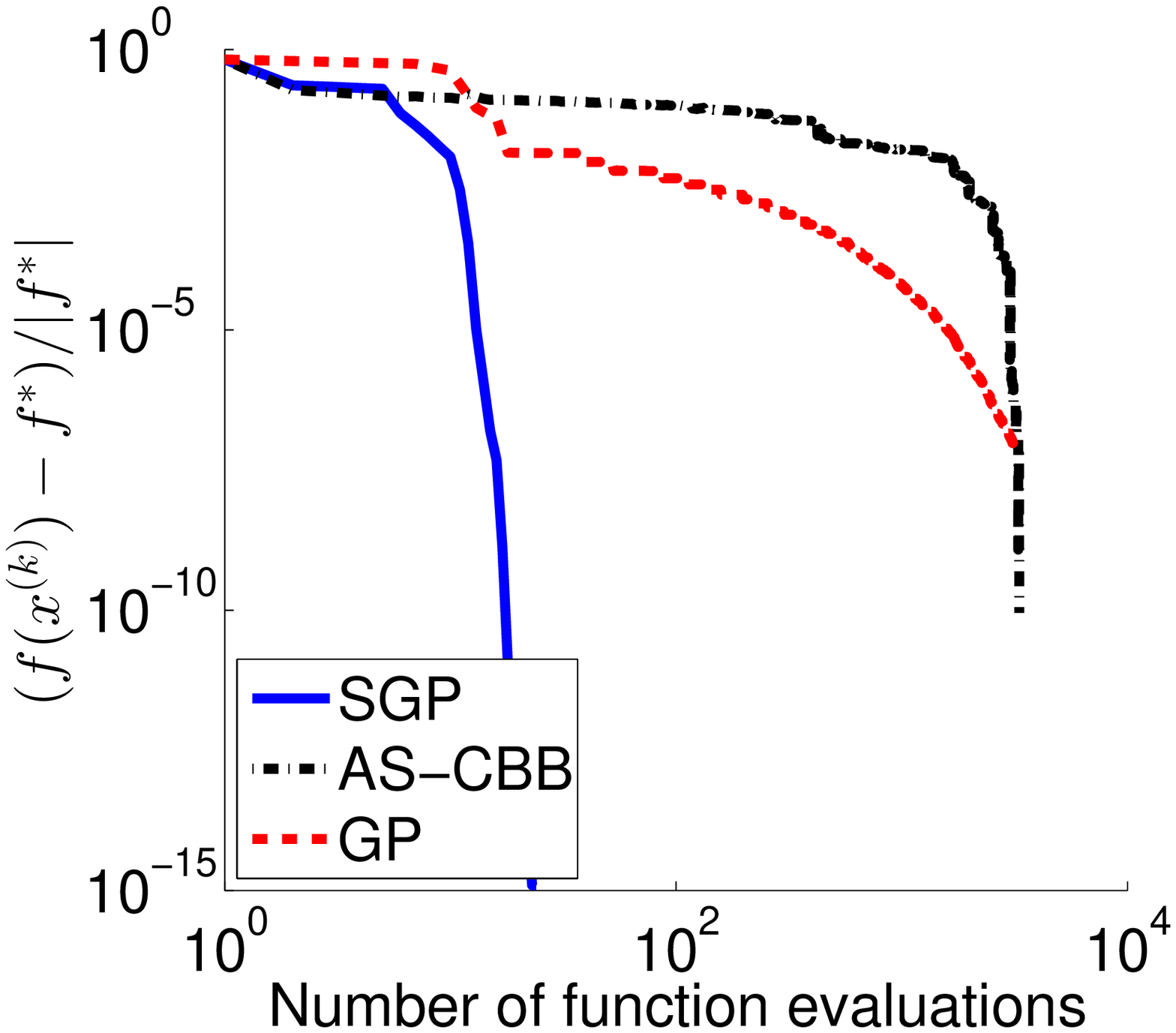}\\
  \includegraphics[width=0.3\linewidth]{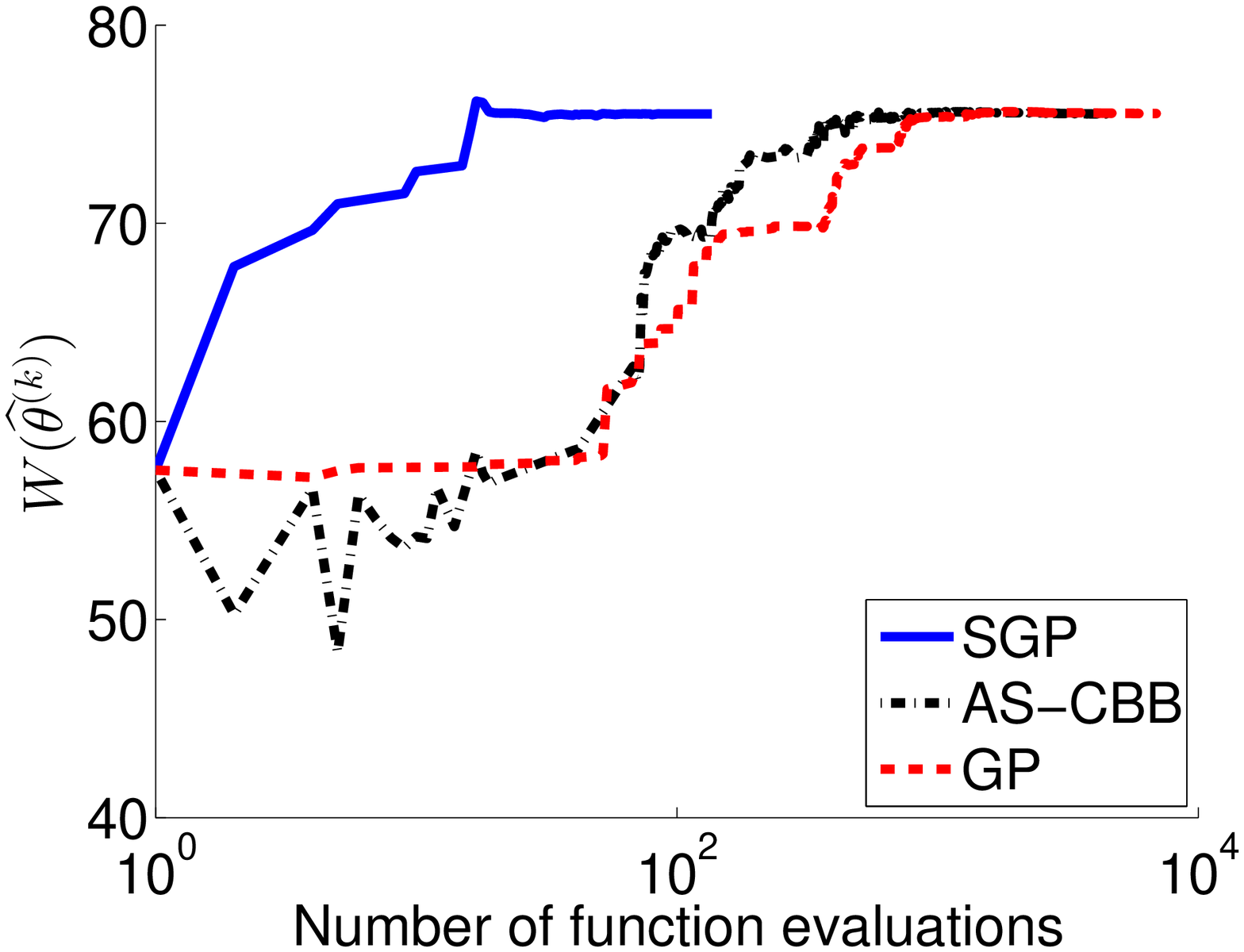} &
  \includegraphics[width=0.3\linewidth]{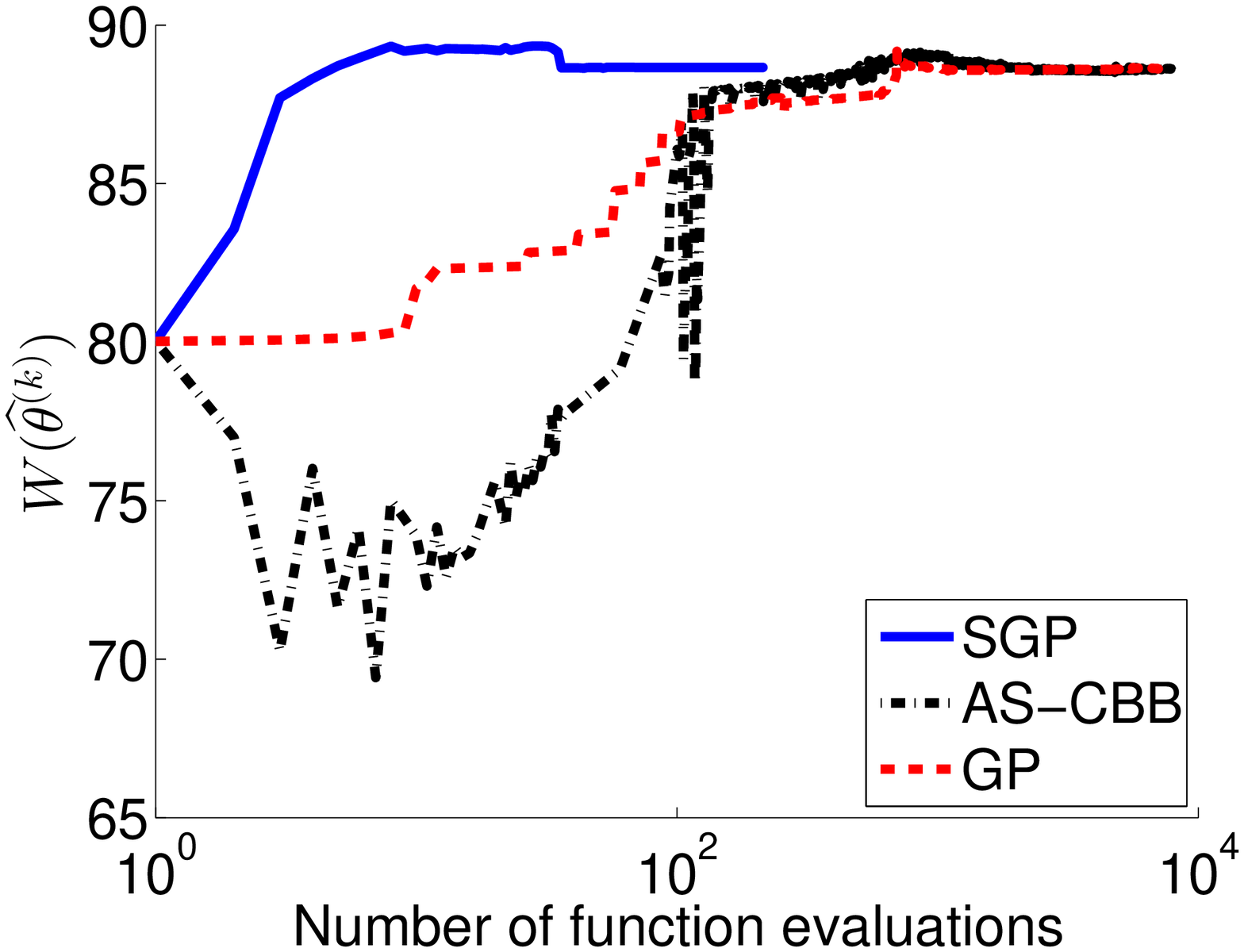} &
	\includegraphics[width=0.3\linewidth]{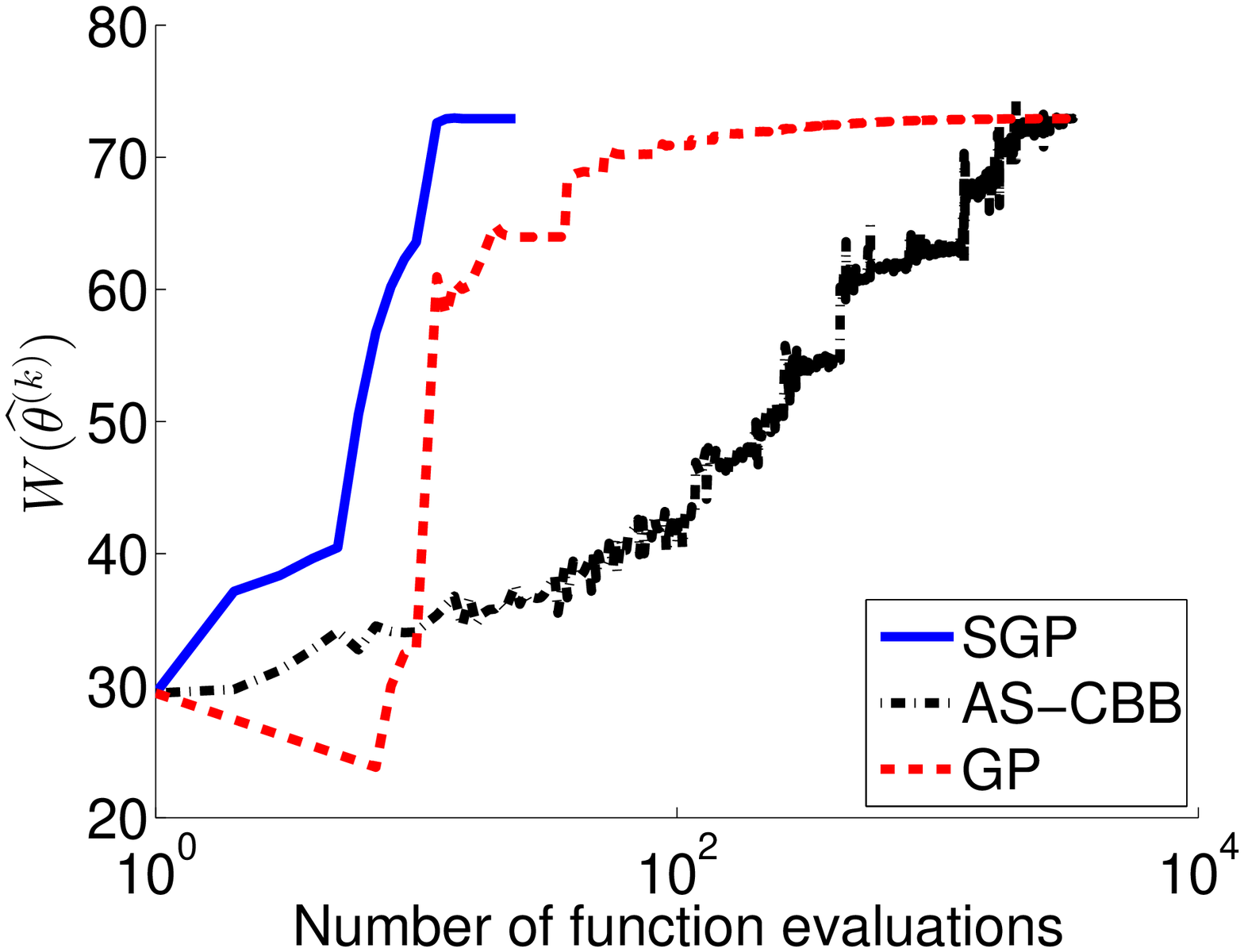}
  \end{tabular}
 \caption{Comparison of SGP, GP and AS-CBB with respect to the number of function evaluations on three instances of the test problems (left: multiple kernel DC-M, dataset D3; middle: multiple kernel TCSS-M, dataset D1; right: kernel SS, dataset D2). First row: relative difference from the reference minimum value. Second row: fit parameter \eqref{fitdef}.}\label{fig:SGP-GP}
\end{center}
\end{figure}

In order to show the significant influence of the scaling strategy in the convergence behaviour of gradient methods, we compare Algorithm \ref{GPM} with the scaling proposed in Section \ref{sec:scaling} (SGP) with the same algorithm without scaling (GP, $D_k = I$) on some instances of the whole test sets described above. Both algorithms adopt the same adaptive alternation of the Barzilai--Borwein rules \eqref{alphaBB} described in Section \ref{subsec:BB} and have all the other parameters set as described in Section \ref{subsec:SGPparams}.\\
As further benchmark, we consider also the affine scaling ciclic Barzilai--Borwein method (AS-CBB) proposed in \cite{Hageretal09}, which consists in a diagonally scaled gradient method whose iteration is given by $\xkk = \xk+\lambda_k\ve d^{(k)}$ where
$$ d^{(k)}_i = \frac{1}{{{\alpha}_k} + |\nabla_i f(\xk)|/X_i(\xk)}\nabla_i f(\xk), \ \ \ \ X_i(\x) = \left\{\begin{array}{ll} u_i -x_i & \mbox{ if } \nabla_i f(\x) \leq 0\\ l_i-x_i & \mbox{ otherwise }\end{array}\right.$$
with the convention $0\cdot \infty = 1/\infty = 0$. In particular, $\alpha_{c\ell+i} = \max(\alpha_{min},1/{\alpha^{BB1}_{c\ell +1}})$, $i=1,...,L_c$ for some fixed cycle length parameter $L_c$ and $\lambda_k$ is computed by a nonmonotone Armijo-type backtracking procedure. When $f$ is twice continuously differentiable, any limit point of the sequence generated by AS-CBB is a stationary point \cite[Theorem 4.1]{Hageretal09}; moreover, the authors also shows the local R-linear convergence to non degenerate local minimum satisfying the second order optimality conditions \cite[Theorem 7.1]{Hageretal09}. In our experiments, we set $L_c = 4$ and the nonmonotone Armijo parameter ($M$ in formula (2.2) in \cite{Hageretal09}) equal to $8$.\\
The plots in Figure \ref{fig:SGP-GP} are obtained by: a) running the Matlab function \verb"fmincon" to get a reference value $f^*$ for a minimum of $f$; b) running each algorithm and computing the relative difference between $f^*$ and the current estimate $f(\xk)$ at each iterate.\\
A significantly faster decrease of the objective function value is observed for SGP, with respect to the number of function evaluations, together with a smoother and faster improvement of the estimated impulse response, measured by means of the fit parameter defined in \eqref{fitdef}. In practice, after the very first SGP iterations, a good estimate of the impulse response is obtained.\\
The comparison between SGP, GP and AS-CBB  gives information about the relative behaviour of scaled gradient methods (SGP, AS-CBB) with respect to a non scaled one (GP) and also about the importance of the scaling matrix choice (SGP versus AS-CBB), which, as observed before, leads to very different performances.

\subsection{SGP results and performance assessment}

\begin{figure}[ht]
\begin{center}
\begin{tabular}{ccc}
  \includegraphics[width=0.3\linewidth]{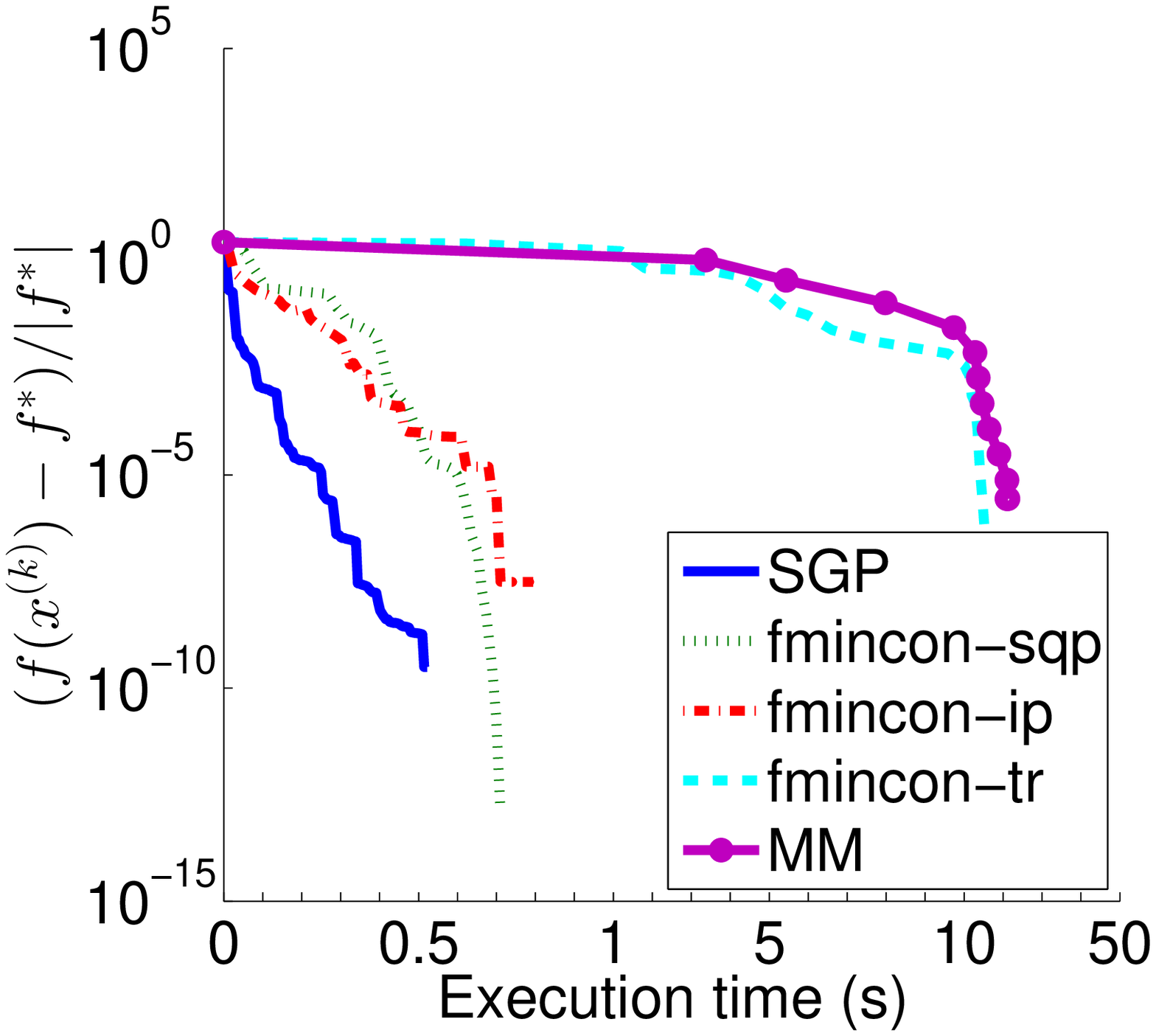} &
  \includegraphics[width=0.3\linewidth]{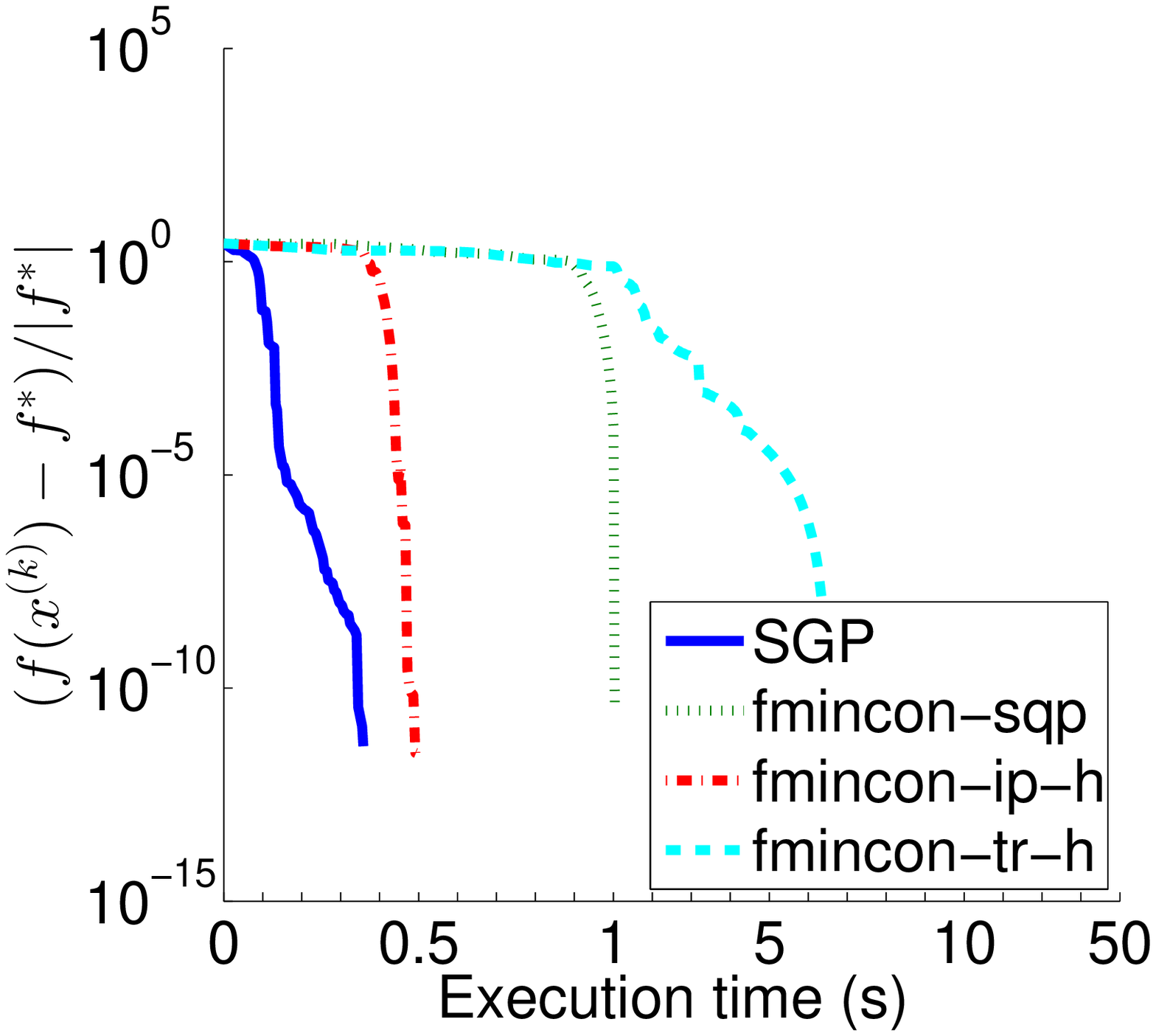} &
	\includegraphics[width=0.3\linewidth]{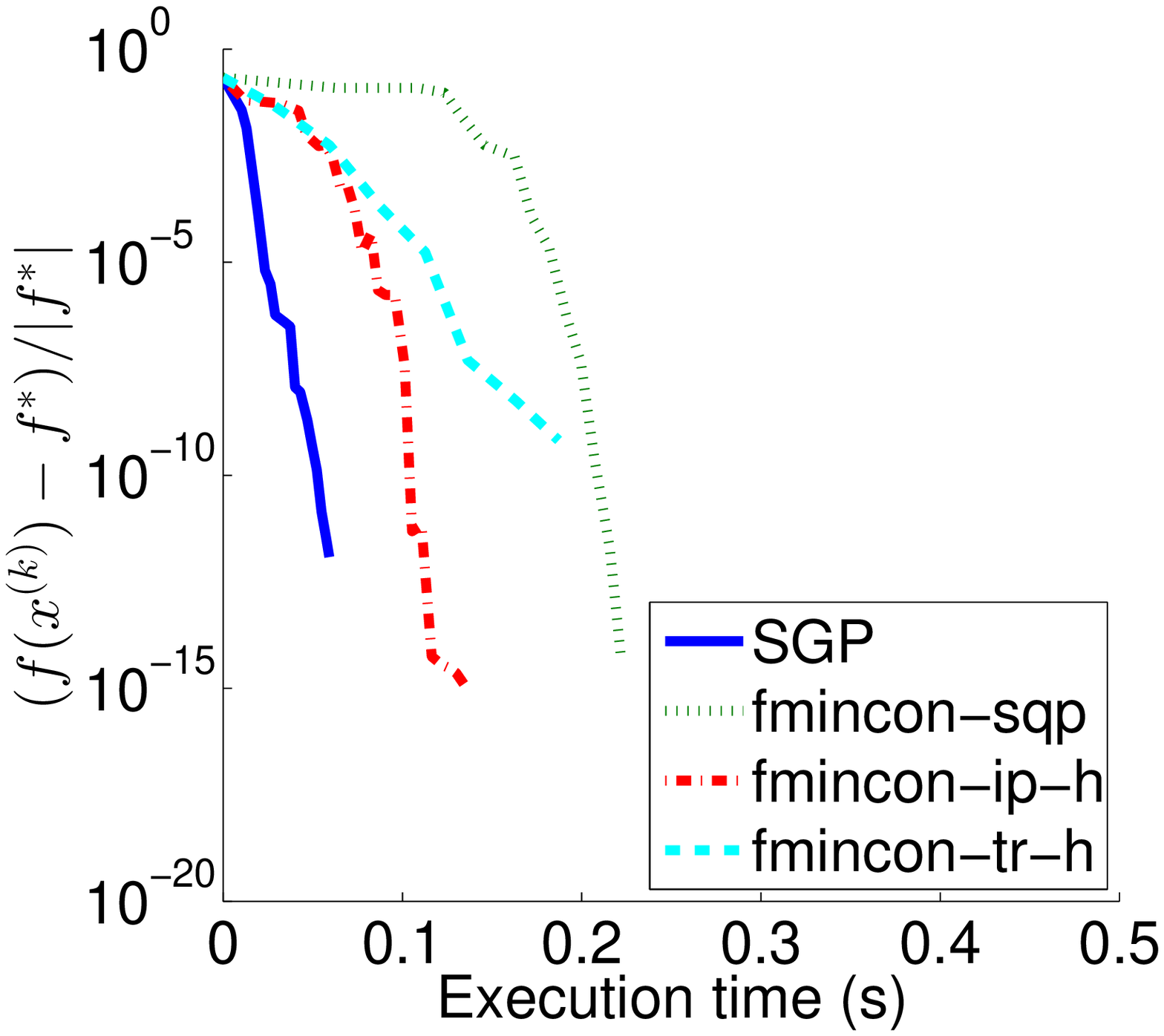}\\
  \includegraphics[width=0.3\linewidth]{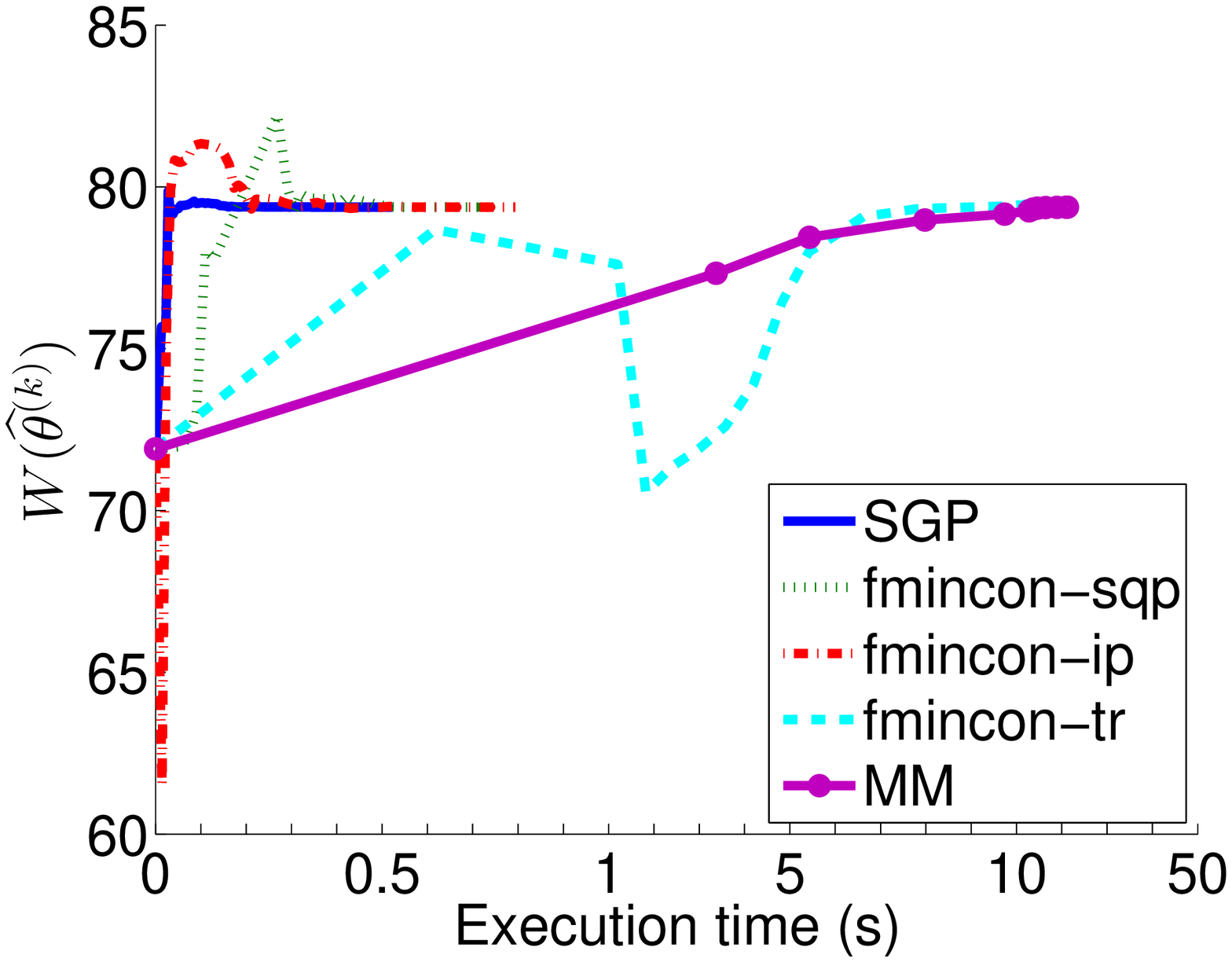} &
  \includegraphics[width=0.3\linewidth]{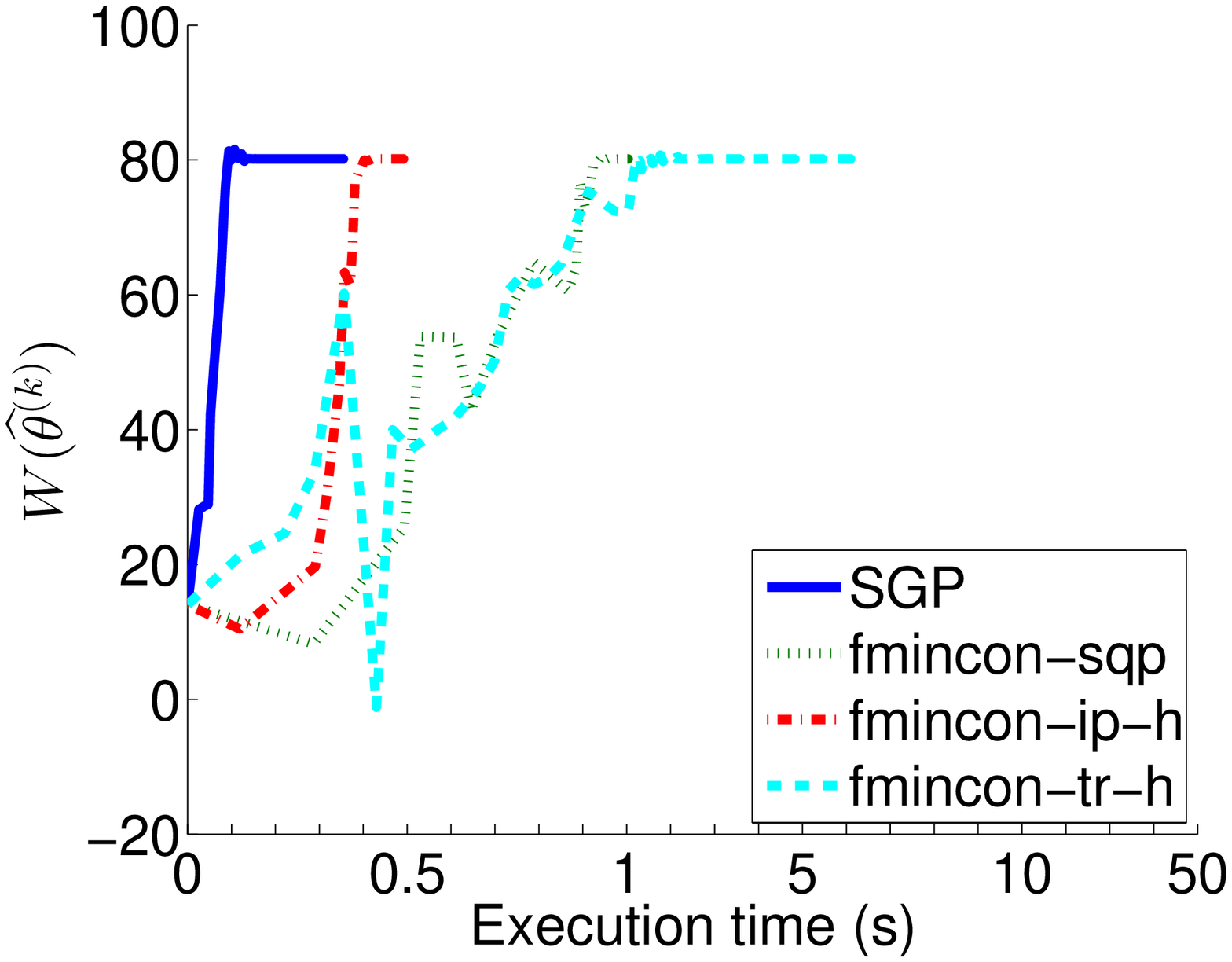} &
	\includegraphics[width=0.3\linewidth]{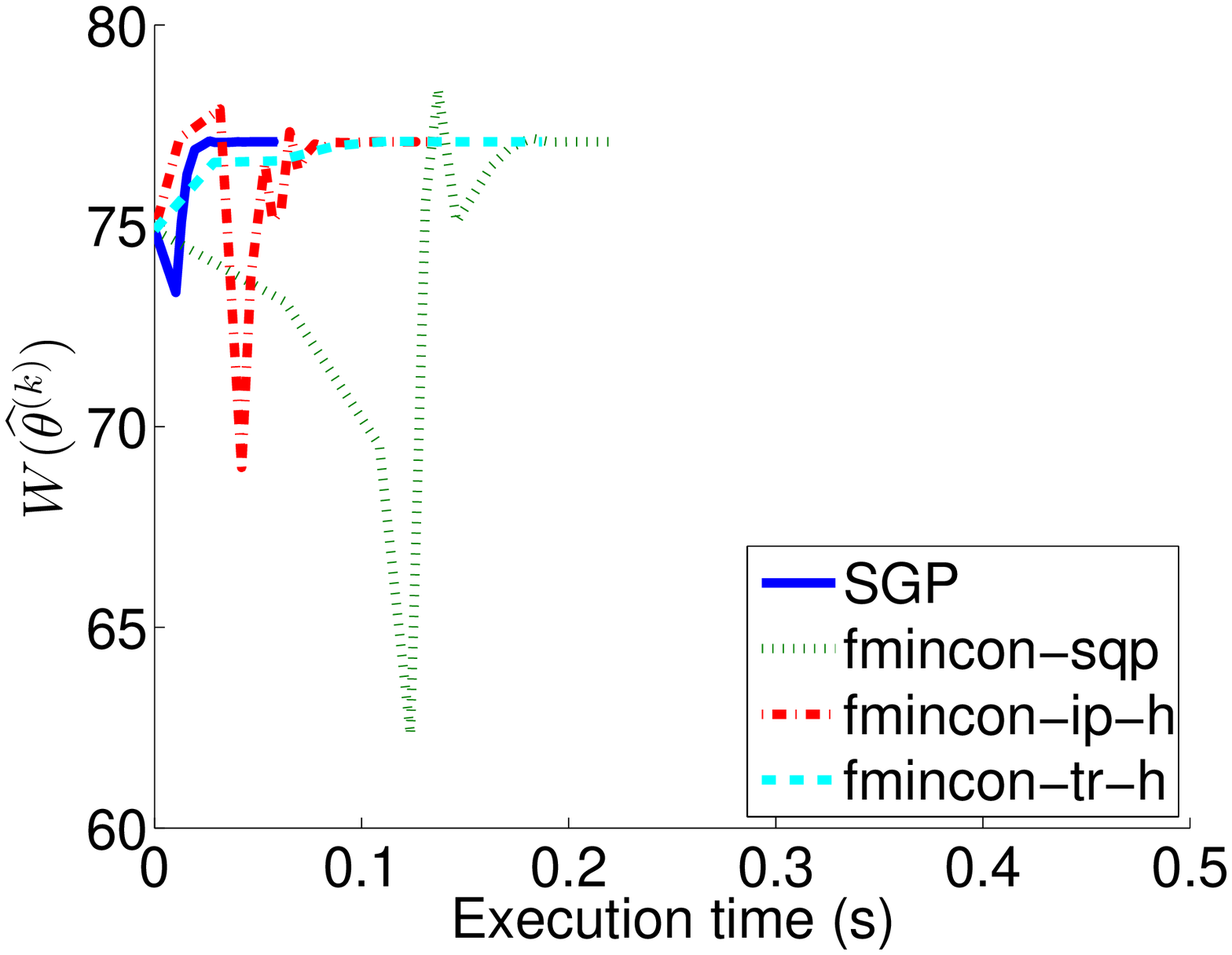}
  \end{tabular}
 \caption{Comparison of SGP and fmincon (with three different algorithm options) with respect to the execution time on three instances of the test problems (left: multiple kernel DC-M, dataset D1; middle: kernel SS, dataset D1; right: kernel TC, dataset D4). The MM algorithm is also shown in the multiple-kernel case. First row: relative difference from the reference minimum value. Second row: fit parameter \eqref{fitdef}.}\label{fig:SGP-GP2}
\end{center}
\end{figure}

\begin{figure}
\begin{center}
\begin{tabular}{cc}
\includegraphics[width=0.3\linewidth]{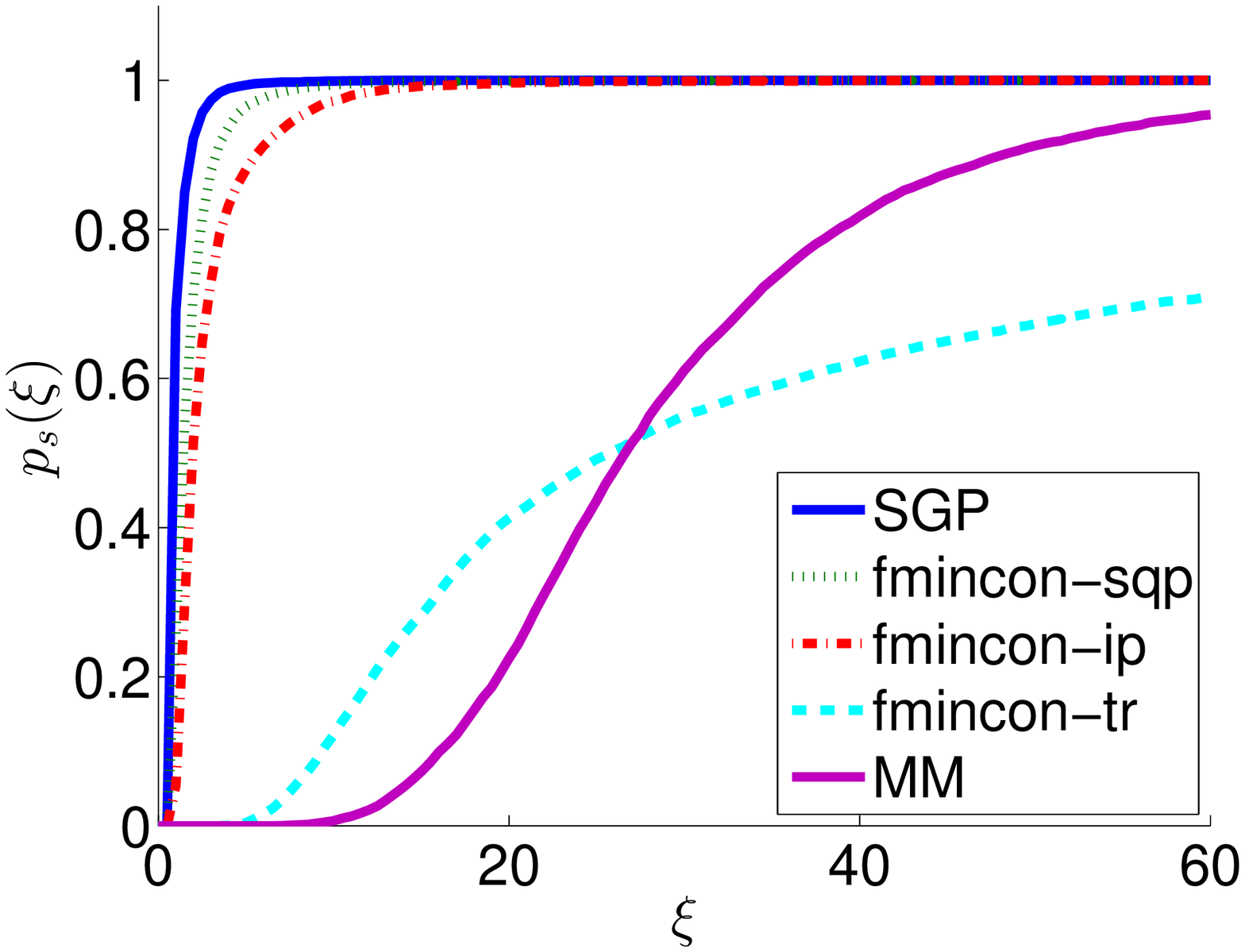} &
\includegraphics[width=0.3\linewidth]{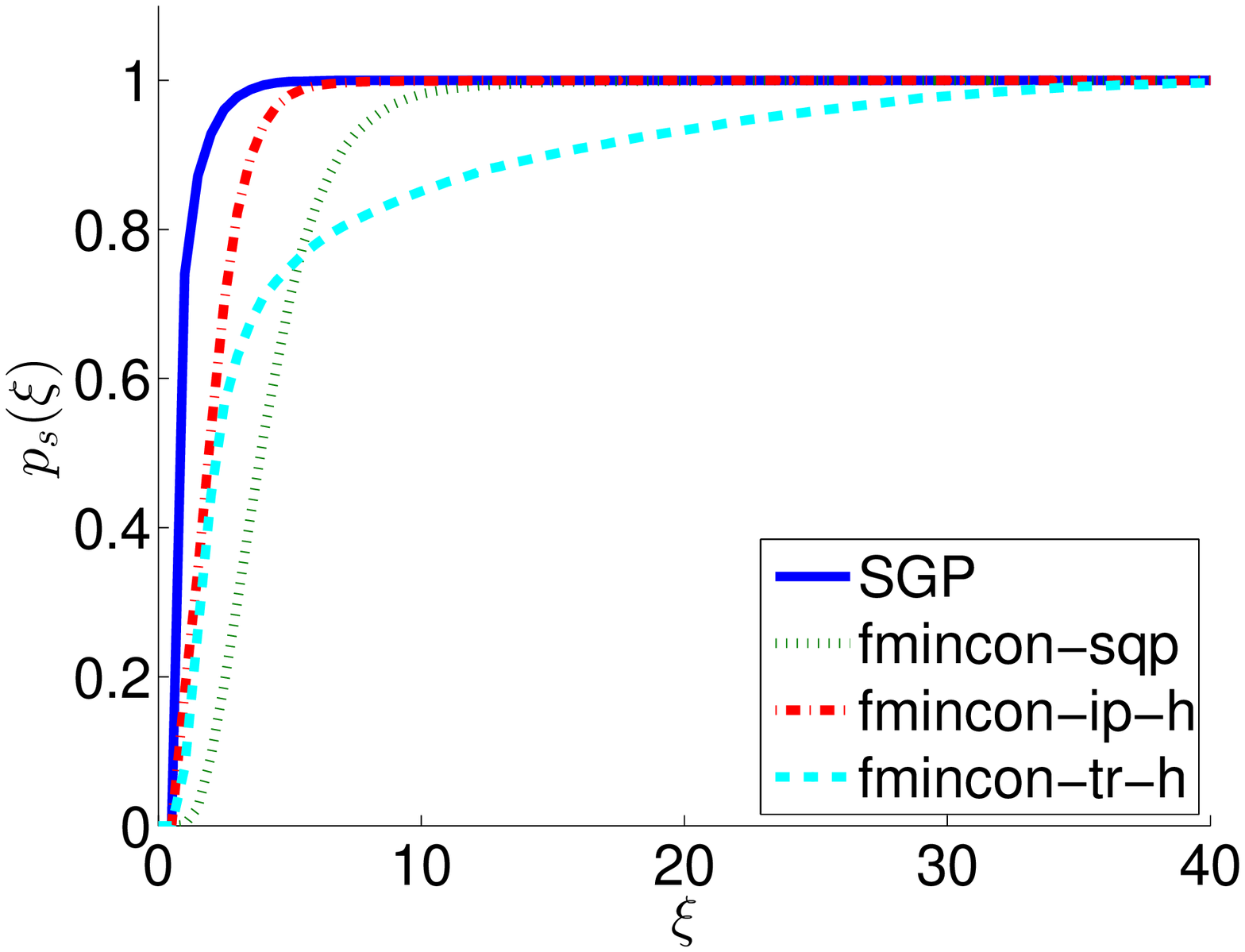}
\end{tabular}
\caption{Performance profiles. Left: multiple kernels DC-M and TCSS-M. Right: kernels DC, TC and SS.}\label{fig:performance-profile}
\end{center}
\end{figure}

In Tables \ref{tab:linear} and \ref{tab:nonlinear} we summarize the results obtained by applying SGP on the test sets described above. For each dataset, we report the average fit \eqref{fitdef}, the average number of iterations (`it'), the average number of function evaluations (`nf') and the average computational time in seconds (`t'). The whole experimentation has been performed with the Matlab implementation of SGP described in the previous section, running on a server with a dual Intel Xeon QuadCore E5620 processor at 2,40 GHz, 12 Mb cache and 18 Gb of RAM under Matlab2010b. The accuracy of the results in terms of the fit parameter \eqref{fitdef} is coherent with the results reported in \cite{ChenetalTAC:14}. To evaluate the effectiveness of SGP, we compare it to other state-of-the-art methods, such as the optimization algorithms implemented in the Matlab \verb"fmincon" function \verb"sqp", \verb"interior-point" and \verb"trust-region-reflective", which is the default one, denoted by `sqp', `ip' and `tr' respectively in the tables. We point out that \verb"fmincon" has a further algorithm option, \verb"active-set", which is however not suited for the considered problems since it may produce infeasible iterates outside the objective function domain.\\
The optimization parameters for \verb"fmincon" are the default ones except \verb"TolFun" which has been set to $10^{-9}$, while the same functions described in Section \ref{sec:implementation} and employed also by SGP have been exploited for objective function and gradient evaluations.\\
The \verb"sqp" option correspond to a BFGS approximation of the Hessian matrix, while \verb"interior-point" and \verb"trust-region-reflective" admit also a user supplied Hessian instead of an automatically computed approximation of it. With the suffix `-h' we indicate that the exact Hessian was also provided to the solvers. For sake of brevity, in Tables \ref{tab:linear} and \ref{tab:nonlinear} we only report the case with the best average computational time. Indeed, as observed at the end of Section \ref{sec:implementation}, the Hessian computation is quite costly in the multiple kernel case, so that the time needed for computing it is not balanced by the reduction of the iteration number that one expects when uses the exact second order information. For example, with the DC-M kernel on the dataset D1, the average iterations and function evaluations numbers were 23 and 25 respectively for the \verb"interior-point" with the exact Hessian, but the corresponding average computational time was 11.04 seconds. Some instances of the plots of the relative difference from the minimum value and the fit parameter \eqref{fitdef} as functions of the execution time for the different strategies are shown in Figure \ref{fig:SGP-GP2}. For the multiple kernel case we also consider the method recently proposed in \cite{ChenetalTAC:14}, which is based on a Minimization-Majorization (MM) approach. In practice, this method solves a sequence of convex optimization subproblems whose objective function is obtained by linearizing the concave term. Each subproblem, which can be formulated as a second order cone program, is then solved by an especially tailored interior point method. The theoretical convergence properties of the MM method (see \cite[Theorem 4]{Sriperumbudur-Lankriet-2009}) are substantially identical to that stated in Theorem \ref{teorema}: every limit point of the sequence is stationary. For the multiple kernel \eqref{Palpha}, property \eqref{coercive} of the objective function only guarantees the existence of limit points. We adopt the MM Matlab implementation provided by the authors, which exploits the CVXOPT package, a Python module for convex optimization \cite{cvxopt}. The numerical comparison with MM has been carried out with Python 2.7.1 installed and with the ATLAS library compiled and optimized for our architecture. All methods are initialized with the same starting point.\\
In order to give a more intuitive insight of the comparison among the different solvers, in Figure \ref{fig:performance-profile} we also report the performance profiles \cite{Dolan-More2002} obtained by grouping the test problems according to the kernel type, multiple or single. Given a test set $\mathcal P$ and a set of solvers $\mathcal S$, let us denote by $t_{p,s}$ the computational time required by solver $s\in \mathcal S$ to solve the problem $p\in\mathcal P$. Then, the performance ratio is defined as $\rho_{p,s} = t_{p,s}/\min\{t_{p,s},s\in \mathcal S\}$. When a solver $s$ does not succeed on a problem $p$, the corresponding ratio $t_{p,s}$ is set to a value $\rho_{max} $ such that $\rho_{p,s} \leq \rho_{max}$ for all $p\in \mathcal P$ and $s\in \mathcal S$ and $\rho_{p,s}=\rho_{max}$ if and only if a failure occurred. The performance profile of the solver $s\in \mathcal S$ is $p_{s}(\xi) = \mbox{size}\{p\in{\mathcal P}: \rho_{p,s} \leq \xi\}/\mbox{size}\{\mathcal P\}$, for a given $\xi\in\R$. The quantity $p_{s}(\xi)$ expresses the probability that a performance ratio $\rho_{p,s}$ lies within a factor of $\xi$ of the best possible ratio.
\begin{table}\footnotesize
\begin{center}
\setlength{\tabcolsep}{5pt}
\begin{tabular}{lc| r   r  rrr| r  r  rrr}
                   &      &\multicolumn{5}{c|}{DC-M}  & \multicolumn{5}{c}{TCSS-M}\\\hline
                   &      &{SGP    }  &MM      &\multicolumn{3}{c|}{fmincon} &{SGP    }  &MM   &\multicolumn{3}{c}{fmincon}\\
                   &      &          &       & sqp     & ip       & tr     &           &        &sqp    & ip        & tr\\\hline
\hline
\multirow{4}{*}{D1}&fit   &      84.4&  84.4 &     84.4&     84.4 &    84.4&      84.4 &84.4    &   84.4&     84.4  &    84.4\\
                   &it    &       104&  12   &       44&       78 &     323&       123 &12      &     42&       72  &      97\\
                   &nf    &       137&  -    &      121&       95 &     324&       156 &-       &    110&       88  &      98\\
                   &t     &      0.74&  16.5 &     0.83&     1.02 &  112.45&      0.58 &9.13    &   0.51&     0.60  &   12.68\\
\hline
\multirow{4}{*}{D2}&fit   &      63.2&  63.1 &     63.2&     63.1 &    63.1&      63.6 &63.6    &   63.6&     63.6  &    63.6\\
                   &it    &        76&  11   &       49&      114 &     126&        94 &11      &     50&       99  &      87\\
                   &nf    &       103&  -    &      121&      134 &     127&       129 &-       &    116&      118  &      88\\
                   &t     &      0.55&  14.6 &     0.83&     1.35 &   44.00&      0.46 &7.89    &   0.54&     0.79  &   11.41\\
\hline
\multirow{4}{*}{D3}&fit   &      87.6&  87.6 &     87.5&     85.7 &    87.6&      88.7 &88.8    &   88.8&     88.6  &    88.8\\
                   &it    &        86&  11   &       52&       99 &     179&       127 &11      &     52&       90  &     129\\
                   &nf    &       115&  -    &      131&      119 &     180&       163 &-       &    131&      108  &     130\\
                   &t     &      0.62&  16.5 &     0.90&     1.25 &   62.48&      0.60 &8.80    &   0.60&     0.73  &   16.89\\
\hline
\multirow{4}{*}{D4}&fit   &      74.9&  74.8 &     74.8&     71.5 &    74.8&      76.6 &76.7    &   76.6&     76.5  &    76.6\\
                   &it    &        76&  10   &       70&      180 &     131&        95 &10      &     68&      139  &      97\\
                   &nf    &       103&  -    &      160&      204 &     132&       126 &-       &    156&      160  &      98\\
                   &t     &      0.55&  14.3 &     1.09&     2.00 &   45.58&      0.46 &7.96    &   0.72&     1.06  &   12.70\\
\hline
\end{tabular}
\caption{Results obtained by SGP, fmincon (with three different algorithm options - see text for details) and MM on multiple kernels DC-M and TCSS-M. For each dataset, we report the average fit \eqref{fitdef}, the average number of iterations (`it'), the average number of function evaluations (`nf') and the average computational time in seconds (`t').}
\label{tab:linear}
\end{center}
\end{table}

\begin{table}\footnotesize
\begin{center}
\setlength{\tabcolsep}{5pt}
\begin{tabular}{ l     c|        r      rrr|    r  rrr| r rrr}
                   &      &\multicolumn{4}{c|}{DC}&\multicolumn{4}{|c|}{TC}&\multicolumn{4}{|c}{SS}\\\hline
                   &      &{SGP}  &  \multicolumn{3}{c|}{fmincon}& SGP &\multicolumn{3}{c|}{fmincon}&SGP &\multicolumn{3}{c}{fmincon}\\
                   &      &       &  sqp   & ip-h   & tr-h   &        &  sqp   & ip-h&tr-h&       &  sqp   & ip-h   &tr-h\\
\hline
\multirow{4}{*}{D1}&fit   &  83.2 &    83.2&    83.2&    83.2&   82.5 &   82.5 & 82.5&82.5&   77.6&    77.5&    77.1&    76.2\\
                   &it    &   124 &      27&      21&     136&     19 &     18 &   13&  17&     43&      30&      21&     225\\
                   &nf    &   168 &      88&      30&     137&     23 &     61 &   22&  18&     62&      81&      31&     226\\
                   &t     &  0.44 &    0.29&    0.25&    1.03&   0.06 &   0.19 & 0.14&0.11&   0.16&    0.27&    0.22&    1.36\\
\hline
\multirow{4}{*}{D2}&fit   &  60.3 &    58.0&    60.2&    60.1&   60.4 &   60.0 & 60.6&60.5&   52.2&    49.8&    52.7&    50.9\\
                   &it    &    59 &      30&      19&      42&     20 &     24 &   15&  16&     29&      29&      21&      96\\
                   &nf    &    77 &      86&      26&      43&     23 &     70 &   22&  17&     38&      75&      28&      97\\
                   &t     &  0.21 &    0.28&    0.22&    0.34&   0.06 &   0.21 & 0.14&0.10&   0.10&    0.25&    0.20&    0.59\\
\hline
\multirow{4}{*}{D3}&fit   &  87.9 &    87.8&    87.8&    87.8&   87.6 &   87.6 & 87.6&87.6&   87.2&    87.1&    87.1&    86.6\\
                   &it    &   109 &      29&      23&     141&     21 &     19 &   14&  18&     40&      33&      24&     262\\
                   &nf    &   148 &      90&      30&     142&     25 &     63 &   24&  19&     55&      86&      34&     263\\
                   &t     &  0.39 &    0.30&    0.25&    1.08&   0.07 &   0.20 & 0.14&0.11&   0.15&    0.29&    0.23&    1.58\\
\hline
\multirow{4}{*}{D4}&fit   &  74.7 &    74.5&    74.8&    74.7&   74.7 &   74.6 & 74.7&74.7&   71.7&    70.5&    72.2&    70.2\\
                   &it    &    78 &      37&      25&      66&     19 &     27 &   17&  16&     32&      34&      23&     172\\
                   &nf    &   104 &     101&      30&      67&     23 &     75 &   26&  17&     41&      83&      31&     173\\
                   &t     &  0.27 &    0.35&    0.29&    0.52&   0.06 &   0.24 & 0.16&0.10&   0.11&    0.28&    0.22&    1.04\\
\hline
\end{tabular}
\caption{Results obtained by SGP and fmincon (with three different algorithm options - see text for details) on DC, TC, SS kernel matrices. For each dataset, we report the average fit \eqref{fitdef}, the average number of iterations (`it'), the average number of function evaluations (`nf') and the average computational time in seconds (`t').}
\label{tab:nonlinear}
\end{center}
\end{table}

From Tables \ref{tab:linear} and \ref{tab:nonlinear} and Figure \ref{fig:performance-profile} we can observe what follows:
\begin{itemize}
\item in general, all the considered methods provide solutions with comparable accuracy, measured in terms of the fit parameter \eqref{fitdef}. Some differences in accuracy could be due to the fact that problem \eqref{pr1} is nonconvex and, then, different algorithms can be attracted by different stationary points; however, the overall results are satisfactory;
\item in presence of simple constraints, a first order method as SGP, equipped with a suitable combination of a scaling matrix and a steplength parameter, is competitive with more sophisticated and highly optimized second order methods, as the ones implemented in the \verb"fmincon" Matlab function;
\item the high flexibility of SGP allows to overcome some limits of state-of-the-art schemes as the MM approach and be applied also when the objective function is not a difference of convex functions.
\end{itemize}

\section{Concluding remarks}\label{sec:conc}

In this paper we have considered linear system identification in the Bayesian framework. A key step is the estimation of the hyperparameters describing the Bayesian prior, which leads to the nonconvex, nonlinear, bound constrained optimization problem \eqref{pr1}. Our aim was to analyze problem \eqref{pr1} from a numerical point of view, proposing also an especially tailored version of SGP for its solution and presenting the results of an extensive numerical experimentation comparing several state-of-the-art algorithms. Our analysis, together with the experimental results, aims to give new insights about the numerical issues related to the considered application and also about gradient projection methods and related scaling techniques. The numerical results, depicted in Figure \ref{fig:performance-profile}, show that the proposed method obtains the overall best performances in terms of time.\\
The main contributions of the paper are summarized below.
From the optimization point of view:
\begin{itemize}
\item we proposed a new split gradient approach for bound constrained optimization;
\item the numerical experience shows that scaling techniques are useful to improve the performances of the gradient projection method also on nonconvex problems. The improvements obtained with the proposed approach are observed with respect to the nonscaled version of the same method and also with respect to a gradient method based on a different scaling technique;
\item the combination of the proposed scaling technique with a suitable steplength selection rule makes SGP competitive with second-order methods, as the ones implemented in the \verb"fmincon" Matlab function.
\end{itemize}
From the application point of view:
\begin{itemize}
\item we provide an ${\mathcal O}(n^3)$ algorithm to evaluate the objective function and gradient of problem \eqref{pr1} and an ${\mathcal O}(mn^3)$ algorithm for Hessian computation;
\item we also provide an extensive numerical experimentation with the Matlab \verb"fmincon" function, devising the most convenient algorithm options.
\end{itemize}
As concluding remarks, we point out that one of the main strength of the proposed approach is the capability to provide a good estimate of the impulse response coefficients after very few iterations, without need of the second order information, which, especially in the multiple kernel case \eqref{Palpha}, is quite costly to compute. We believe that the good performances of SGP rely on the fact that the proposed scaling technique takes into account of the problem structure, that is both the objective function and the constraints. On the other side, this is also the main difficulty to the generalization of SGP: indeed, the scaling technique is especially tailored for box constraints and the extension to more general constraints is not straightforward. This issue will be addressed in our future work, which will consider also a wider range of problems arising in machine learning and system identification where sparse Bayesian learning ideas can be applied, e.g. the identification of multi-input, multi-output systems, where also automatic variable selection needs to be performed, or the basis selection problem in the context of machine learning.


\end{document}